\documentclass[11pt]{article}
\usepackage[margin=1.25in]{geometry}
\usepackage{setspace}
\usepackage{amsmath,amssymb,amsthm,mathtools}
\usepackage{graphicx}
\usepackage{natbib}
\usepackage{enumitem}
\usepackage{algorithm}
\usepackage{algorithmic}
\usepackage{booktabs}
\usepackage[hidelinks]{hyperref}
\usepackage{float}

\newcommand{\A}{\mathcal{A}}
\newcommand{\Asa}{\mathcal{A}_{\mathrm{sa}}}
\newcommand{\R}{\mathbb{R}}
\newcommand{\dd}{\mathrm{d}}

\newcommand{\Wt}{W_{2}}
\newcommand{\pv}{\mathrm{p.v.}\!}
\newcommand{\Ent}{\chi}
\newcommand{\Free}{\mathcal{F}}
\newcommand{\Div}{\mathcal{D}}
\newcommand{\relF}{I}
\newcommand{\norm}[1]{\lVert #1\rVert}
\newcommand{\abs}[1]{\lvert #1\rvert}
\newcommand{\Econd}[1]{\mathcal{E}_{#1}}
\DeclareMathOperator{\supp}{supp}

\theoremstyle{plain}
\newtheorem{theorem}{Theorem}[section]
\newtheorem{proposition}[theorem]{Proposition}
\newtheorem*{proposition*}{Proposition}

\newtheorem{corollary}[theorem]{Corollary}
\theoremstyle{definition}
\newtheorem{example}[theorem]{Example}
\newtheorem{algorithmenv}[theorem]{Algorithm}
\theoremstyle{remark}
\newtheorem{remark}[theorem]{Remark}

\begin{document}
\sloppy

\title{Beyond the Semicircle: Free Diffusion Models with Prescribed Equilibria}
\author{Swagatam Das\\[2pt]
\small Electronics and Communication Sciences Unit, Indian Statistical Institute, Kolkata, India\\
\small\texttt{swagatam.das@isical.ac.in}}
\date{}
\maketitle

\begin{abstract}
A growing class of machine-learning objects -- covariance and Gram matrices,
kernel and attention matrices, MIMO channel matrices, density operators -- are
naturally spectra rather than coordinate vectors. Building a denoising diffusion
model for such data by corrupting eigenvalues coordinatewise is not merely
inelegant: it converges to the wrong limit, because eigenvalues repel rather than
move independently. Free probability theory supplies the correct forward process,
with free convolution replacing classical convolution and Voiculescu's conjugate
variable replacing the score, but constant-coefficient free diffusions have a hidden
limitation of their own: their only possible equilibrium is the semicircular law,
whatever the target distribution looks like. We show that state-dependent free
volatility removes this restriction. For any sufficiently regular compactly
supported target law, we explicitly construct, through a closed-form
Hilbert-transform drift, a free Fokker--Planck flow having that law as a stationary
spectral distribution. Under an additional convex-potential condition, the
construction admits a globally relaxing diffusion interpretation, reverse-time
dynamics, and a trainable score. The generalization is genuine rather than
cosmetic: no change of spectral variable reduces it to the constant-coefficient
case, and it is a Wasserstein gradient flow only for an explicitly characterized
family of coefficients that contains no bounded non-constant member. We instantiate
the construction on the Marchenko--Pastur law, the limiting spectrum of isotropic
sample covariance matrices, and validate it numerically: the free-versus-coordinatewise
gap, an exact dimension-independent score driving reverse-time matrix dynamics at
several matrix sizes, the designed Marchenko--Pastur equilibrium, and end-to-end
generation with a learned score.
\end{abstract}

\section{Introduction}\label{sec:intro}

Denoising diffusion models generate samples by learning to reverse a stochastic
process that gradually destroys structure \citep{SohlDickstein2015, Ho2020,
Song2021}. The theory is built for vector-valued data: a forward
Ornstein--Uhlenbeck flow corrupts $x_0\in\R^d$ with Gaussian noise, and a
reverse-time SDE driven by the score $\nabla\log p_t$ regenerates it.

A growing share of the data machine learning actually operates on is not naturally a
coordinate vector at all, but a \emph{spectrum}: the eigenvalues of a sample
covariance or Gram matrix, a kernel or attention matrix, a MIMO channel matrix, or a
quantum density operator. For such objects the natural invariance is unitary
conjugation, not coordinate permutation, and the natural summary statistics --
capacity, generalization bounds, entropy, purity -- are functionals of the spectrum
alone. The standard machinery of diffusion models -- forward SDE, score, reverse-time
equation -- is not stated in a form that respects this invariance. This raises a basic
question: what does it mean to build a diffusion model whose target is a spectrum, and
does the vector-valued theory even give the right answer if applied naively?

\textbf{It does not.} The obvious approach is to treat the $N$ eigenvalues as $N$
exchangeable scalars and corrupt them with independent Gaussian noise, exactly as
in the vector-valued theory. For unitarily invariant spectral data this is provably the wrong limit. Eigenvalues of a
Hermitian matrix perturbed by Hermitian noise do not move independently -- they
repel -- and the correct large-$N$ limit is governed by \emph{free convolution}
\citep{VDN1992}, not classical convolution. The two make quantitatively different
predictions: in our experiments the coordinatewise approximation is off by two
orders of magnitude relative to the free one (\S\ref{sec:experiments}).

\begin{figure}[t]
\centering
\includegraphics[width=0.92\textwidth]{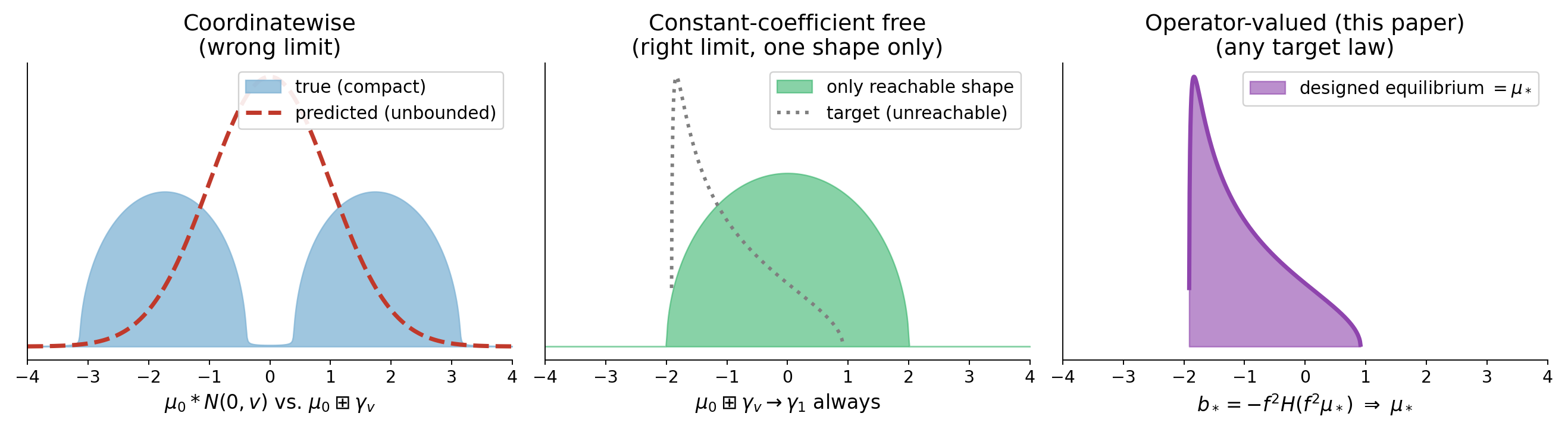}
\caption{Three ways to corrupt a spectrum, and what each can produce as an
equilibrium. Left: coordinatewise noise gives the wrong limit entirely. Middle: the
correct (free) process, but only ever semicircular. Right: this paper -- an
arbitrary designed target.}
\label{fig:warmup}
\end{figure}

Free probability theory, Voiculescu's noncommutative analogue of classical
probability, supplies the correct forward process, its own Fokker--Planck equation,
its own notion of score (the \emph{conjugate variable}), and a reverse-time SDE
\citep{Dabrowski2014, BianeSpeicher1998}. But building this theory on the direct
analogue of the vector-valued forward process -- a free Ornstein--Uhlenbeck flow with
scalar (state-independent) diffusion coefficient -- inherits a second, less obvious
restriction: \emph{the only possible equilibrium of such a flow is the semicircular
law}, up to rescaling. Whatever the target spectrum actually is -- Marchenko--Pastur
for a sample covariance matrix, or whatever an attention or channel matrix produces
-- the constant-coefficient construction can only ever generate the semicircle.
Figure~\ref{fig:warmup} illustrates all three regimes at a glance.

\paragraph{Scope.} The object generated throughout is a \emph{spectral law}, the
eigenvalue distribution of a unitarily invariant matrix ensemble. Eigenvector information
is intentionally outside the model: this is a generative model of spectral distributions,
not of arbitrary matrices.

\paragraph{Contributions.} Three main results, following Figure~\ref{fig:warmup}, and a
group of structural results.

\emph{1. The correct spectral diffusion.} For unitarily invariant spectral data,
coordinatewise corruption yields the classical limit $\mu_0*\mathcal N(0,v)$ and
Hermitian corruption the free limit $\mu_0\boxplus\gamma_v$; they differ at the fourth
moment by exactly $v(2m_2+v)$ and by two orders of magnitude in $W_1$ at $N=2500$
(\S\ref{sec:naive}).

\emph{2. Equilibrium design (the central result).} Letting the diffusion coefficient
depend on the state removes the semicircular restriction: for any sufficiently regular
compactly supported target $\mu_*$ and positive weight $f$ (Theorem~\ref{thm:design}),
\[
  \boxed{\;b_*(x)=-f(x)^2\,H(f^2\mu_*)(x)\;}
\]
makes $\mu_*$ a stationary spectral law of
$\dd X_t=b_*(X_t)\,\dd t+f(X_t)\,\dd S_t\,f(X_t)$. We keep three levels of claim apart
(\S\ref{sec:design}): stationarity of the continuity equation, for every regular target;
relaxation, reversal and a trainable score, for strictly convex potentials
(Corollary~\ref{cor:designconvex}; Marchenko--Pastur is an instance); and well-posedness
of the SDE, not established in general.

\emph{3. A generative realization.} A free Tweedie identity (\S\ref{sec:background})
reduces score estimation to regression of a scalar function of eigenvalue and noise
level, whose parameterization is independent of the matrix size $N$;
\S\ref{sec:experiments} validates it with an exact score at several $N$ and a learned
score at one $N$.

\emph{Structural results.} No diffeomorphism of the spectral line reduces non-constant
state-dependent volatility to a constant one (Theorem~\ref{thm:lamperti}); the flow is a
Wasserstein gradient flow only for an explicitly characterized family of coefficients
containing no bounded non-constant member (Theorem~\ref{thm:nogradient}); and an
explicit threshold governs when a spectral gap survives partway through the forward
corruption (Theorem~\ref{thm:gapclosing}). Proofs and further experiments are in the
appendix.

\section{Background}\label{sec:background}

This section fixes notation for \S\ref{sec:design}; readers who want the main
construction first can skip ahead after Figure~\ref{fig:warmup} and return to the
free-entropy material only when Corollary~\ref{cor:designconvex} is reached. We work in a tracial
$W^*$-probability space $(\A,\tau)$: a von Neumann algebra $\A$
with a faithful normal tracial state $\tau$, playing the role of expectation. A
self-adjoint $X\in\Asa$ has a \emph{spectral law} $\mu_X$, the unique probability
measure on $\R$ with $\tau(X^n)=\int x^n\,\dd\mu_X(x)$ for all $n$; this is the
noncommutative analogue of the law of a random variable, and it is the only object
our diffusion models act on. \emph{Free independence} replaces classical
independence, and \emph{free convolution} $\boxplus$ replaces classical convolution
as the operation governing sums of freely independent variables. A \emph{free
Brownian motion} $(S_t)$ is the free analogue of Brownian motion: $S_t$ has the
semicircular law $\gamma_t$ with density $\tfrac{1}{2\pi t}\sqrt{4t-x^2}$. See
\citet{VDN1992, NicaSpeicher2006} for the standard references.

\paragraph{Forward process.} The free Ornstein--Uhlenbeck SDE
\begin{equation}\label{eq:ffwd}
  \dd X_t = -\tfrac12\beta(t)X_t\,\dd t + \sqrt{\beta(t)}\,\dd S_t
\end{equation}
is the free analogue of the classical forward diffusion. Its spectral marginals solve
a \emph{nonlocal} Fokker--Planck equation carried by the Hilbert transform,
\begin{equation}\label{eq:ffp}
  \partial_t\psi_t = \tfrac{\beta}{2}\partial_x(x\psi_t) - \beta\,\partial_x(\psi_t
  H\mu_t),
\end{equation}
where $H\mu_t(x)=\pv\!\int\frac{\dd\mu_t(y)}{x-y}$ is the Hilbert transform,
a complex Burgers equation, solved explicitly and uniquely, in the class of Cauchy
transforms of probability measures (this is the right uniqueness class for the
problem, not merely a convenient one; see Appendix~\ref{app:diperna}), by
\begin{equation}
\mu_t=D_{\sqrt{\alpha_t}}\mu_0\boxplus\gamma_{1-\alpha_t}
\end{equation}
with $\alpha_t=e^{-\int_0^t\beta}$; for every $t>0$ this $\mu_t$ has a bounded,
real-analytic density on the interior of its support (Biane's regularity theorem),
regardless of the regularity of $\mu_0$. As $t\to\infty$,
$\mu_t\to\gamma_1$, the semicircular law, regardless of $\mu_0$: this is the
restriction our main result removes.

\paragraph{Free score and reverse-time SDE.} Voiculescu's \emph{conjugate variable}
$\xi_{\mu_t}:=2H\mu_t$ is the free analogue of the negative score $-\nabla\log p_t$
(both equal $x$ for the standard Gaussian, respectively semicircular, law);
Table~\ref{tab:dictionary} gives the dictionary. Specializing the time reversal of free diffusions
\citep{Dabrowski2014} to \eqref{eq:ffwd} gives the reverse-time free SDE
\begin{multline}\label{eq:reversesde}
  \dd Y_s = \big[\tfrac12\beta(T\!-\!s)Y_s \\
  - \beta(T\!-\!s)\xi_{\mu_{T-s}}(Y_s)\big]
  \dd s + \sqrt{\beta(T\!-\!s)}\,\dd\bar S_s,
\end{multline}
which transports the semicircular law back to $\mu_0$ (a second, independent
derivation via finite-$N$ Haussmann--Pardoux reversal gives the same equation in the
$N\to\infty$ limit; the two routes and exactly how they relate are compared in
Appendix~\ref{app:tworoutes}). A free Tweedie identity
expresses $\xi_{\mu_t}$ as a conditional expectation, $\xi_{\mu_t}(Y)=v_t^{-1}(Y-
\Econd{W^*(Y)}[\sqrt{\alpha_t}X_0])$, turning score estimation into a regression
problem with an exact consistency guarantee (deferred to
Appendix~\ref{app:tweedie}), the free analogue of denoising score matching
\citep{Vincent2011, Hyvarinen2005}.

\begin{table}[t]
\centering
\caption{Dictionary between classical and free spectral diffusion.}
\label{tab:dictionary}
\small
\begin{tabular}{@{}p{0.42\columnwidth}p{0.52\columnwidth}@{}}
\toprule
Classical (vectors) & Free (spectra)\\
\midrule
Gaussian convolution & free convolution $\boxplus\gamma_v$\\
score $-\nabla\log p_t$ & conjugate variable $\xi_{\mu_t}=2H\mu_t$\\
Tweedie, denoising score matching & free Tweedie identity (App.~\ref{app:tweedie})\\
reverse SDE & reverse free SDE \eqref{eq:reversesde}\\
constant-coefficient limit: Gaussian & constant-coefficient limit: semicircle only\\
Langevin drift $\tfrac12\nabla\log p_*$ & $b_*=-H\mu_*=-\tfrac12\xi_{\mu_*}$ ($f\equiv1$)\\
\bottomrule
\end{tabular}
\end{table}

\paragraph{Free entropy and functional inequalities.} Two functionals organize the
convexity theory this paper relies on. Voiculescu's free entropy and the free energy
relative to the quadratic potential are
\begin{align}
  \Ent[\mu] &= \iint\log\abs{s-t}\,\dd\mu(s)\dd\mu(t)+\tfrac34+\tfrac12\log2\pi,\\
  \Free[\mu] &= \tfrac12\!\int\! x^2\dd\mu(x) - \Ent[\mu],
\end{align}
and, writing $\Div(\mu\Vert\gamma):=\Free[\mu]-\Free[\gamma]$ for the relative free
energy and $\relF(\mu\Vert\gamma):=\int(x-\xi_\mu(x))^2\dd\mu(x)$ for the relative
free Fisher information, $\Free$ is $1$-displacement-convex along free Wasserstein
geodesics \citep{Biane2003}. This single fact yields the free logarithmic Sobolev inequality
$\Div\le\tfrac12\relF$ \citep{Biane2003}, the Talagrand inequality $\Wt^2\le2\Div$
\citep{HiaiPetzUeda2004}, and the HWI inequality $\Div\le\Wt\sqrt{\relF}-\tfrac12\Wt^2$
\citep{LedouxPopescu2009}, each between $\mu$ and $\gamma$ (HWI implies the other two),
together with the de~Bruijn identity
$\tfrac{\dd}{\dd\Lambda}\Div(\mu_t\Vert\gamma)=-\tfrac12\relF(\mu_t\Vert\gamma)$ and
exponential decay $\Div(\mu_t\Vert\gamma)\le e^{-\Lambda(t)}\Div(\mu_0\Vert\gamma)$,
$\Lambda(t)=\int_0^t\beta$. Corollary~\ref{cor:designconvex} shows this entire
apparatus transfers verbatim to any strictly convex confining potential $V$ in place
of $\tfrac12x^2$, which is what makes the equilibrium-design construction of
\S\ref{sec:design} more than an existence statement.

\paragraph{Matrix realization.} All of this is the $N\to\infty$ limit of a concrete,
simulable object: $N\times N$ Hermitian matrix diffusions $\dd X^N_t = b(X^N_t)\,\dd
t + N^{-1/2}\dd H_t$ (GUE noise $H_t$, a Dyson Brownian motion), whose empirical
spectral distribution converges to $\mu_t$ \citep{Dyson1962, RogersShi1993}. This is
what we simulate throughout \S\ref{sec:experiments}.

\section{Why the naive constructions fail}\label{sec:naive}

Two direct ways of building a diffusion model on spectra are available before
introducing anything free-probabilistic. Both fail, for different reasons.

\paragraph{Coordinatewise corruption is wrong, not just inefficient.} Treat the $N$
eigenvalues $\lambda=(\lambda_1,\dots,\lambda_N)$ as exchangeable scalars and corrupt
them with independent Gaussian noise, exactly as the vector-valued theory prescribes.
This gives the classical convolution $\mu_0 * \mathcal N(0,v)$. But eigenvalues of a
Hermitian matrix corrupted by Hermitian noise do not move independently: they repel,
and the correct limit is the free convolution $\mu_0\boxplus\gamma_v$. The two differ
even in gross qualitative features -- $\mu_0\boxplus\gamma_v$ has compact support,
$\mu_0*\mathcal N(0,v)$ does not -- and their fourth moments differ by exactly
$v(2m_2+v)$ (Appendix~\ref{app:freevsclassical}). Table~\ref{tab:freevsclassical}
quantifies the gap: at $N=2500$ the free prediction tracks the truth to
finite-$N$-fluctuation accuracy, while the coordinatewise prediction is off by two
orders of magnitude.

\paragraph{The correct (free) forward process has a hidden restriction.} Once
corruption is modeled correctly by \eqref{eq:ffwd}, the forward flow converges to the
semicircular law $\gamma_1$ for \emph{every} initial law $\mu_0$, because that is the
only equilibrium a constant-coefficient free Ornstein--Uhlenbeck process can have. A
diffusion model built directly on \eqref{eq:ffwd} therefore cannot represent target
spectra such as Marchenko--Pastur, or whatever an attention or channel matrix
produces, except as the byproduct of a schedule that never quite finishes. Making the
equilibrium a design choice is the content of \S\ref{sec:design}.

\section{Designing the equilibrium}\label{sec:design}

We allow the volatility to depend on the current state, using the sandwich
coefficient $f(X)\,\dd S\,f(X)$ -- so that $f(X)$, an operator-valued coefficient built
from $X$ itself, plays the role a scalar diffusion coefficient plays classically
(the unique choice compatible with self-adjointness
of the noise term up to symmetrization; see Appendix~\ref{app:selfadjoint}):
\begin{equation}\label{eq:opsde}
  \dd X_t = b(X_t)\,\dd t + f(X_t)\,\dd S_t\,f(X_t).
\end{equation}

\begin{theorem}[Spectral dynamics under sandwich volatility]\label{thm:opvol}
Let $f$ be bounded and continuous and $b$ bounded and Lipschitz. The spectral law
$\mu_t$ of the solution of \eqref{eq:opsde} satisfies
$\partial_t\mu_t+\partial_x(\mu_t V_t)=0$ with velocity field
\begin{equation}\label{eq:opvelocity}
  V_t(x)=b(x)+f(x)^2\,H\eta_t(x),\qquad \dd\eta_t := f(x)^2\dd\mu_t(x).
\end{equation}
\end{theorem}

\begin{figure}[t]
\centering
\includegraphics[width=0.92\textwidth]{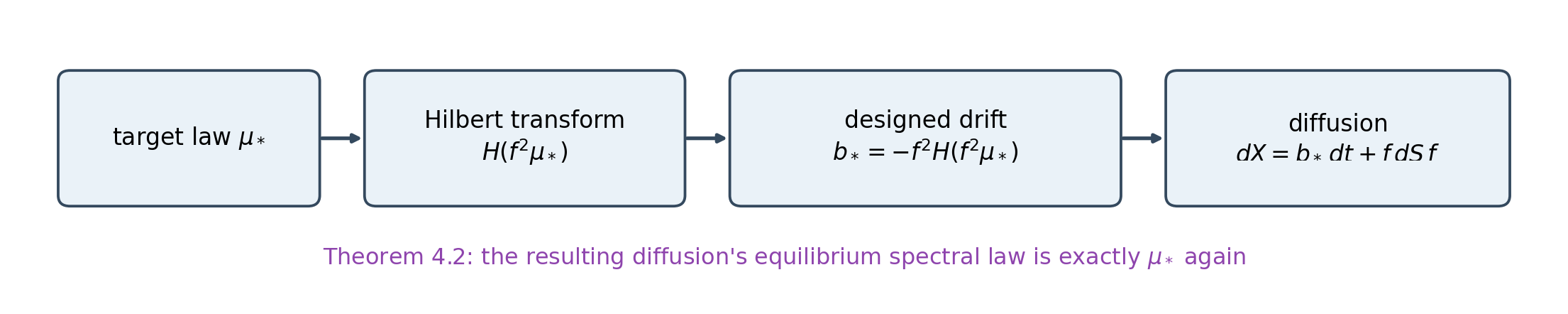}
\caption{The equilibrium-design construction: a target law $\mu_*$ determines,
through a single Hilbert-transform formula, a drift $b_*$ for which $\mu_*$ itself is a
stationary spectral law of the free flow.}
\label{fig:pipeline}
\end{figure}

Constant $f\equiv\beta^{1/4}$ recovers \eqref{eq:ffp}. Setting $V\equiv0$ and solving
for $b$ gives our main result, whose construction Figure~\ref{fig:pipeline}
summarizes end to end.

\paragraph{Three levels of claim.} ``A free diffusion with prescribed equilibrium'' can
mean three things, which we keep apart. \textbf{(I) Stationarity:} the constant path
$\mu_t\equiv\mu_*$ solves the continuity equation of Theorem~\ref{thm:opvol}, i.e.\ $V$
vanishes on $\supp\mu_*$. \textbf{(II) Relaxation:} the flow also relaxes to $\mu_*$
exponentially fast (the gradient-flow regime of Corollary~\ref{cor:designconvex}).
\textbf{(III) Well-posed realization:} the SDE \eqref{eq:opsde} with drift $b_*$ has a
unique solution. Theorem~\ref{thm:design} proves (I) for every regular target,
Corollary~\ref{cor:designconvex} gives (II) for strictly convex potentials, and (III) is
not established in general (\S\ref{sec:limitations}).

\begin{theorem}[Equilibrium design]\label{thm:design}
Let $\mu_*$ be a compactly supported probability measure with a H\"older continuous
density on $K:=\supp\mu_*$, and let $f\in C(K)$, $f>0$. Then
\begin{equation}\label{eq:designdrift}
  b_*(x) = -f(x)^2\,H(f^2\mu_*)(x)
\end{equation}
is well defined and H\"older continuous on compact subsets of the interior of $K$
(and, if $\mu_*$'s density vanishes at the edges of $K$ at a H\"older rate, on all of
$K$), and $\mu_*$ is a stationary solution of the spectral continuity equation of
Theorem~\ref{thm:opvol} for \eqref{eq:opsde} with drift $b_*$ and weight $f$
(Level~I).
\end{theorem}

\begin{proof}[Proof sketch]
By construction $b_*(x)+f(x)^2H(f^2\mu_*)(x)\equiv0$ on the interior of $K$, so the
constant path $\mu_t\equiv\mu_*$ trivially solves the continuity equation of
Theorem~\ref{thm:opvol}. Regularity of $b_*$ follows from the Plemelj--Privalov
theorem for the Hilbert transform of H\"older functions. Full proof, including the
edge case and a companion uniqueness statement, in Appendix~\ref{app:design}.
\end{proof}

\begin{proposition}[Uniqueness of the design drift]\label{prop:designunique}
If in addition $\mu$'s density vanishes at every edge of $K$ (a soft edge, as for
Marchenko--Pastur or the semicircle), then $b_*$ is the \emph{only} bounded continuous
drift, paired with $f$, for which $\mu$ is stationary; and semicircular is the only
stationary law of any constant-coefficient free OU flow, by strict displacement
convexity of $\mathcal F$ (Appendix~\ref{app:inequalities}).
\end{proposition}

\begin{corollary}[Convex potentials]\label{cor:designconvex}
Take $f\equiv1$ so $b_* = -\tfrac12 V'$ with $V'=2H\mu_*$. If $V$ is strictly convex,
$\mu_*$ is the unique minimizer of the strictly displacement-convex free energy
$\mathcal F_V[\mu]=\int V\,\dd\mu - \Ent[\mu]$, the diffusion is its
$\Wt$-gradient flow, and the full relaxation/reverse-time apparatus of
Appendix~\ref{app:inequalities} applies with $V$ in place of $\tfrac12x^2$:
exponential relaxation to $\mu_*$, a reverse-time SDE, and a trainable score.
\end{corollary}

\begin{example}[Marchenko--Pastur]\label{ex:mp}
Let $\mu_*$ be Marchenko--Pastur with ratio $L$, supported on
$[(1-\sqrt L)^2,(1+\sqrt L)^2]$, the limiting spectrum of sample covariance matrices of isotropic data.
Its equilibrium potential $V(x)=\tfrac{x}{L}-\tfrac{1-L}{L}\log x$ satisfies
$V''(x)=\tfrac{1-L}{Lx^2}>0$ on $(0,\infty)$, so Corollary~\ref{cor:designconvex}
applies: $\dd X_t = -\tfrac12V'(X_t)\,\dd t+\dd S_t$ relaxes to Marchenko--Pastur, a
diffusion model whose equilibrium is not semicircular. Well-posedness of the SDE
itself (as opposed to the continuity equation) needs care near $x=0$, where $V'$ is
unbounded; see Appendix~\ref{app:design} for the precise statement and
\citet{WeiYin2026} for the relevant well-posedness theory.
\end{example}

\subsection{Two obstructions: the generalization is genuine}\label{sec:obstructions}

A natural objection is that state-dependent free volatility might secretly reduce to the
constant-coefficient case already understood -- either by a change of spectral
variable, as the classical Lamperti transform does, or by being a Wasserstein
gradient flow after all. Both fail.

\begin{theorem}[No spectral Lamperti reduction]\label{thm:lamperti}
No diffeomorphism $g$ of $\R$ conjugates the spectral flow of \eqref{eq:opsde}, for
every initial law, to a free SDE with constant diffusion coefficient, unless $g$ is
affine and $f$ is constant.
\end{theorem}

\begin{theorem}[Gradient-flow characterization]\label{thm:nogradient}
The spectral flow of \eqref{eq:opsde} is the $m$-mobility Wasserstein gradient flow
of a pairwise-interaction energy for every initial law \emph{if and only if}
\begin{equation}\label{eq:gradflowclass}
  f(x)^4 = \kappa x^2+\lambda x+\nu > 0\ \text{on } \R,\qquad m = f^4/c,\ c>0,
\end{equation}
a two-parameter family of coefficients. Among bounded positive $f$ the only member of
this family is the constant case $\kappa=\lambda=0$, since a positive quadratic is
bounded only when it is constant; the non-constant members, such as
$f(x)=(x^2+1)^{1/4}$, are unbounded.
\end{theorem}

Both proofs (Appendix~\ref{app:obstructions}) trace to the same mechanism: the
Hilbert transform of a pushforward measure is not the pushforward of the Hilbert
transform, so the nonlocal term in \eqref{eq:opvelocity} does not commute with either
reparametrization or the local structure a gradient flow requires. A further
proposition (Appendix~\ref{app:obstructions}) shows that allowing $k$-body
interaction energies for any finite $k$ does not enlarge the family
\eqref{eq:gradflowclass} either. The practical consequence: for the bounded,
non-constant coefficients this paper's well-posedness theory covers
(Theorem~\ref{thm:opvol}), the convexity toolkit of Appendix~\ref{app:inequalities}
is unavailable, and Theorem~\ref{thm:design}'s direct-substitution argument is not a
detour around harder analysis -- it is the only route available.

\subsection{When does a spectral gap survive?}\label{sec:gapclosing}

A third question the design construction raises: if $\mu_*$ has a gap in its support
-- as, e.g., a mixture of well-separated clusters would -- does the forward flow
\eqref{eq:ffwd} that trains a diffusion model on such data close that gap
instantaneously, or does the gap persist for a while? The answer is quantitative.

\begin{theorem}[Gap closing along the free OU flow]\label{thm:gapclosing}
Let $\mu_0$ be compactly supported with a bounded gap $(\alpha_-,\alpha_+)$ in its
support carrying mass on both sides, and let
$\Theta_{\mu_0}(u)=\int(u-y)^{-2}\dd\mu_0(y)$ be its second-order Cauchy kernel. With
$T:=\min_{u\in(\alpha_-,\alpha_+)}\Theta_{\mu_0}(u)\in(0,\infty)$, the corresponding
gap of $\supp\mu_t$ persists exactly while $\alpha_t > T/(1+T)$, closing at
$\Lambda(t)=\log(1+1/T)$.
\end{theorem}

The threshold is a genuine potential-theoretic quantity: $\Theta_{\mu_0}=-G_{\mu_0}'$,
minus the second derivative of $\mu_0$'s logarithmic potential, so a gap with small
curvature (a wide, shallow gap) closes late, and a sharp one closes early
(Appendix~\ref{app:gapclosing}). Proof and a two-independent-gap numerical check are
in Appendix~\ref{app:gapclosing}.

\section{Training and sampling}\label{sec:algorithm}

Equation~\eqref{eq:reversesde} requires $\xi_{\mu_t}$, unavailable in closed form
from data. The free Tweedie identity of \S\ref{sec:background} turns its estimation
into a regression problem: given a denoiser $g_\theta$, define
\begin{equation}\label{eq:loss}
  \mathcal L(\theta)=\int_0^T\!\! w(t)\,\tau\big[(g_\theta(X_t,t)-\sqrt{\alpha_t}X_0)^2
  \big]\dd t,
\end{equation}
whose minimizer recovers $\xi_{\mu_t}$ exactly (Appendix~\ref{app:tweedie}). For
unitarily invariant data the optimal denoiser acts spectrally, so the object learned
is a single scalar function $h_\theta(\lambda,\alpha_t)$ whose parameterization does not
depend on the matrix size $N$, so one score can in principle serve every $N$ (tested
with the exact score in \S\ref{sec:experiments}).

\begin{algorithmenv}[Training and sampling]\label{alg:main}
\textbf{Training.} Sample $X_0^N$, $t\sim\mathrm{Unif}[0,T]$, GUE noise $H$; form
$X_t^N=\sqrt{\alpha_t}X_0^N+\sqrt{v_t}H$; diagonalize $X_t^N=U\Lambda U^*$; regress
the diagonal of $U^*(\sqrt{\alpha_t}X_0^N)U$ on $(\lambda,\alpha_t)$ to fit
$h_\theta$; set $\widehat\xi_t(\lambda)=(\lambda-h_\theta(\lambda,\alpha_t))/v_t$.
\textbf{Sampling.} Draw $X_T^N$ from the GUE; integrate \eqref{eq:reversesde} in its
matrix form backwards by Euler--Maruyama, applying $\widehat\xi_t$ spectrally at
each step and re-symmetrizing.
\end{algorithmenv}

\paragraph{Cost.} The score is dimension-independent in parameterization, not in
computation: each step of the dense implementation diagonalizes an $N\times N$ matrix at
$O(N^3)$ cost, and reducing this is left open.

\section{Experiments}\label{sec:experiments}

All experiments use Hermitian matrix diffusions at the stated $N$; further experimental detail is in
Appendices~\ref{app:experiments} and~\ref{app:realdata}.

\paragraph{Forward dynamics and the hydrodynamic limit.} Starting from the two-atom law
$\pm a_0$ at $N=1200$, the finite-$N$ matrix diffusion matches the free Fokker--Planck
prediction \eqref{eq:ffp} through the support transition below, with $W_1$ error between
$1.2\times10^{-3}$ and $2.4\times10^{-3}$ (Figure~\ref{fig:spacetime},
Appendix~\ref{app:movedfigs}). Averaged over realizations at $N=25,\dots,800$, the error
against the \emph{exact} free marginal decays as $N^{-1}$ (fitted exponents $-0.93$,
$-0.97$), reflecting GUE rigidity; measured against an equal-size Monte~Carlo sample it
shows the misleading rate $N^{-1/2}$, the fluctuation of the reference sample
(Appendix~\ref{app:additional-figures}).

\paragraph{The support transition.}\label{sec:transition-para} For the two-atom law
$\tfrac12(\delta_{-a}+\delta_a)\boxplus\gamma_v$, a bisection search on the numerically
computed density locates the gap-closing transition at $v^*=2.5600$ against the
theoretical $a_0^2=2.56$, and the density at the gap center over the $(a,v)$ plane
confirms $v^*=a^2$ as an identity in $a$ (Appendix~\ref{app:additional-figures}).

\paragraph{Functional inequalities.} Along the same forward trajectory we verify the free
log-Sobolev, Talagrand and HWI inequalities of \S\ref{sec:background} as ratios that must
stay below $1$: the maxima are $0.491$, $0.522$ and $0.678$ (HWI tightest, as it implies
the other two), and the de~Bruijn identity holds to a median relative error of
$2.4\times10^{-3}$ over four orders of magnitude in $\Div$ (Figure~\ref{fig:lsi},
Appendix~\ref{app:movedfigs}).

\paragraph{Gap survival: two independent gaps.} For a three-atom law at
$-2.3,0.4,1.9$ with weights $0.5,0.2,0.3$, Theorem~\ref{thm:gapclosing} predicts two
\emph{different} thresholds, $\alpha_1^*=0.287$ and $\alpha_2^*=0.479$, since $T$
depends on position and mass, not only on gap width. Counting connected components
of the numerically computed support confirms both, in the predicted order and within
one grid step of each threshold (Appendix~\ref{app:gapclosing}) -- direct evidence
that the theory controls not just whether but exactly \emph{when} structure in a
target spectrum survives partway through the forward corruption, which is precisely
the information a practitioner needs to choose how much of the schedule to run.
(This does not cover the vanishing-mass outliers of the spiked model, which follow
Baik--Ben~Arous--P\'ech\'e-type phenomena \citep{BaikBenArousPeche2005}.)

\begin{figure}[h]
\centering
\includegraphics[width=0.85\linewidth]{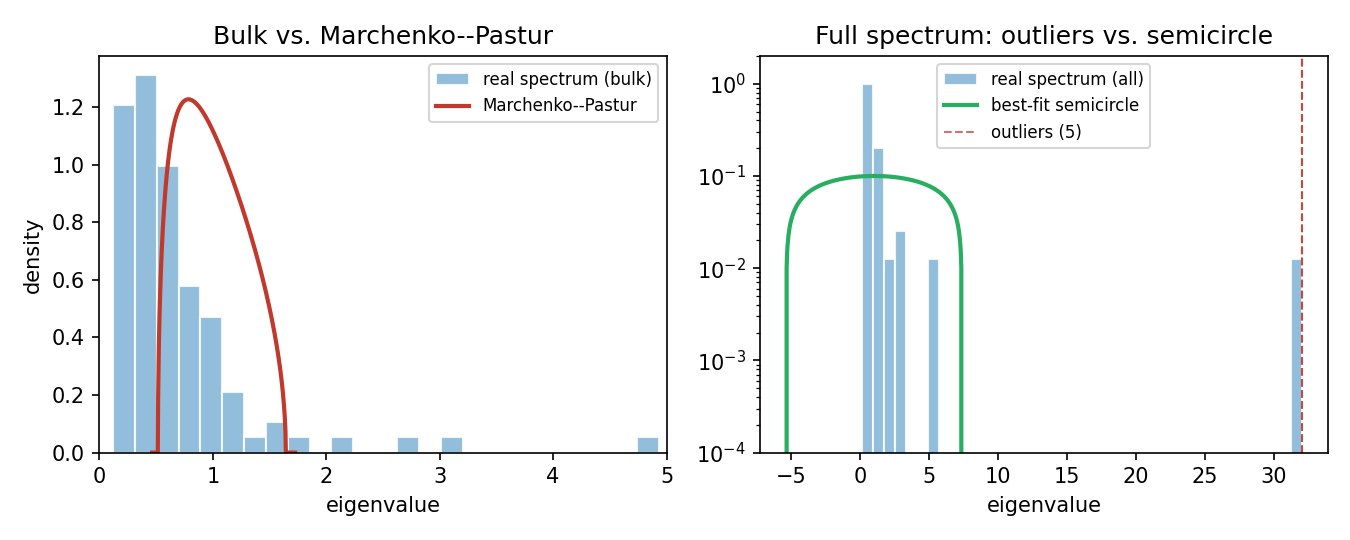}
\caption{Real S\&P~500 correlation spectrum ($N=100$, $T=1258$) against the
Marchenko--Pastur bulk prediction and the best-fitting semicircular law -- the only
shape a constant-coefficient free diffusion can produce.}
\label{fig:realdata}
\end{figure}

\paragraph{A real spectral distribution outside the semicircular family.} We use daily
closing prices of the $N=470$ S\&P~500 constituents with complete coverage from
2013-02-08 to 2018-02-07 \citep{plotlydatasets}, take standardized log-returns, and form
the empirical correlation matrix of $N=100$ assets over all $T=1258$ trading days
($L=N/T=0.080$; selection and preprocessing in Appendix~\ref{app:realdata}). As is
familiar for financial correlations \citep{Laloux1999,Plerou2002,Bun2017}, the spectrum
is only loosely Marchenko--Pastur-like -- the $95$ smallest eigenvalues span $[0.12,1.68]$
against a predicted support $[0.52,1.64]$, and only $47\%$ lie inside it -- and $5$
eigenvalues sit clearly above the predicted edge, the largest carrying $32\%$ of the
total variance (the market mode). No constant-coefficient free diffusion can represent
this: its only equilibrium is semicircular, which even at the best-fitting mean and
variance misses the real spectrum by $W_1=2.20$ (Figure~\ref{fig:realdata}).
Theorem~\ref{thm:design} applies not to the raw eigenvalues but to a regularized density
estimate $\hat\mu_*$, a truncated, renormalized kernel density estimate of the eigenvalues
(Appendix~\ref{app:realdata}): $b_*=-H\hat\mu_*$ makes $\hat\mu_*$, not the raw empirical
measure, a stationary spectral law (Level~I). Simulating the designed diffusion from a
mismatched start moves the bulk toward the target region, but fully resolving the
isolated market-mode eigenvalue is slow (eigenvalue repulsion) and left as future work.
This experiment shows that real spectra can lie far outside the constant-coefficient
equilibrium family; it is not a demonstrated generative application.

\paragraph{Free versus coordinatewise corruption.} We corrupt a two-atom spectrum
($\pm1.6$) with Hermitian Gaussian noise of variance $v$ at $N=2500$ and compare the
empirical spectrum against both predictions.

\begin{table}[h]
\centering
\caption{$W_1$ distance to the empirical spectrum. Free tracks the finite-$N$
fluctuation; coordinatewise is off by two orders of magnitude.}
\label{tab:freevsclassical}
\begin{tabular}{lccc}
\toprule
$v$ & $0.25$ & $1.00$ & $2.56$ \\
\midrule
free ($\boxplus$) & 0.0010 & 0.0013 & 0.0022 \\
coordinatewise ($*$) & 0.109 & 0.219 & 0.275 \\
\bottomrule
\end{tabular}
\end{table}

\begin{figure}[t]
\centering
\includegraphics[width=0.85\textwidth]{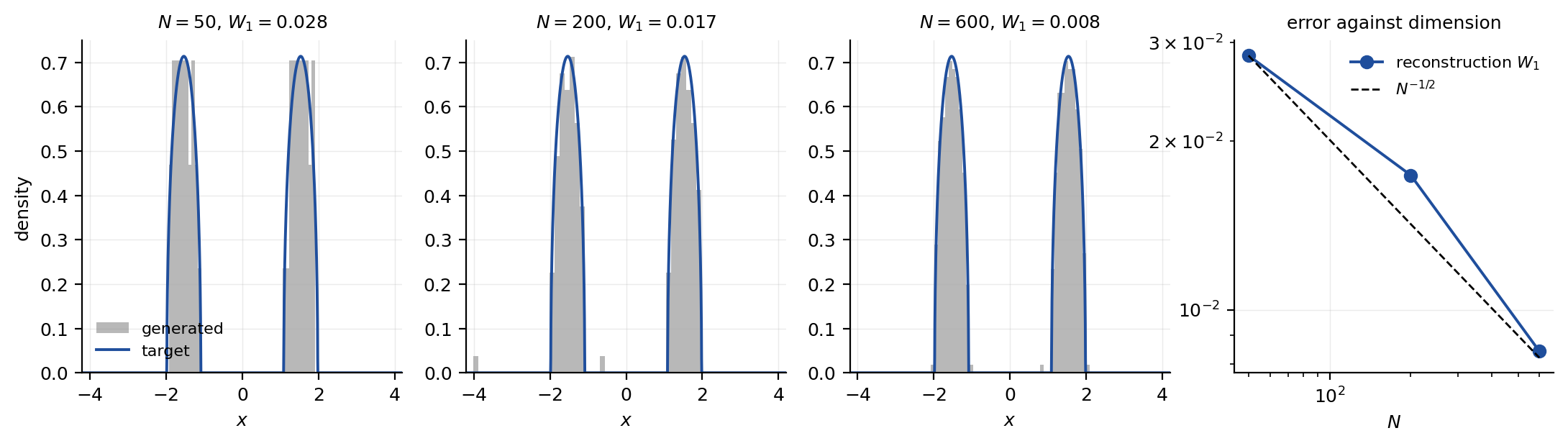}
\caption{Reverse-time matrix dynamics driven by the exact free score of the two-atom
law at $N=50,200,600$ against the target density, and $W_1$ error against $N$.}
\label{fig:transfer}
\end{figure}

\paragraph{One score for every dimension.} The free score is a scalar function of the
eigenvalue and the noise level, so the same function can drive the reverse-time matrix
SDE at any $N$. To test this in isolation from estimation error we use the \emph{exact}
free score of the two-atom law, computed once from its closed-form density, at
$N=50,200,600$: the reconstruction error is $W_1=0.028,\,0.017,\,0.0084$ (fitted
exponent $-0.48$, consistent with the $O(N^{-1/2})$ accumulation of Euler--Maruyama
noise; Figure~\ref{fig:transfer}). This does not yet test a \emph{learned} score across
sizes: the learned score below is evaluated at a single $N$. An entrywise score, by
contrast, is an $N^2$-dimensional object tied to its training size.

\paragraph{A designed equilibrium: Marchenko--Pastur.} We integrate the Hermitian
matrix diffusion of Example~\ref{ex:mp} at $N=400$, $L=\tfrac12$, from an initial
spectrum far from the target. The eigenvalue system itself is stiff near the origin
(the Dyson repulsion grows like the reciprocal gap while $V'$ diverges there), so we
integrate the Hermitian matrix model instead, applying $V'$ spectrally at each step so
repulsion is generated implicitly by the matrix noise (Theorem~\ref{thm:opvol}) rather
than formed explicitly, stable at a step size where the direct eigenvalue integration
is not (Appendix~\ref{app:experiments}). The
empirical spectrum converges to Marchenko--Pastur. We report two checks that bear
directly on the hypotheses of Theorem~\ref{thm:design}, rather than a single
terminal-distance summary: the stationarity identity $V'=2H\mu_*$ holds to relative error
$3.5\times10^{-3}$ on the interior of the support, and the fitted edge exponents of
the empirical density are $0.42$ and $0.51$ against the theoretical $1/2$, confirming
the H\"older edge-vanishing that makes $b_*$ bounded on all of $K$.

\paragraph{Generative modeling with a learned score.} We instantiate
Algorithm~\ref{alg:main} end to end on a spiked-covariance spectrum at $N=400$,
comparing three scores driving the same sampler.

\begin{table}[h]
\centering
\caption{$W_1$(generated, data), same reverse-time sampler, three scores.}
\label{tab:generative}
\begin{tabular}{lccc}
\toprule
score & learned & oracle & coordinatewise \\
\midrule
$W_1$ & 0.073 & 0.062 & 0.121 \\
\bottomrule
\end{tabular}
\end{table}

\begin{figure}[h]
\centering
\includegraphics[width=0.85\linewidth]{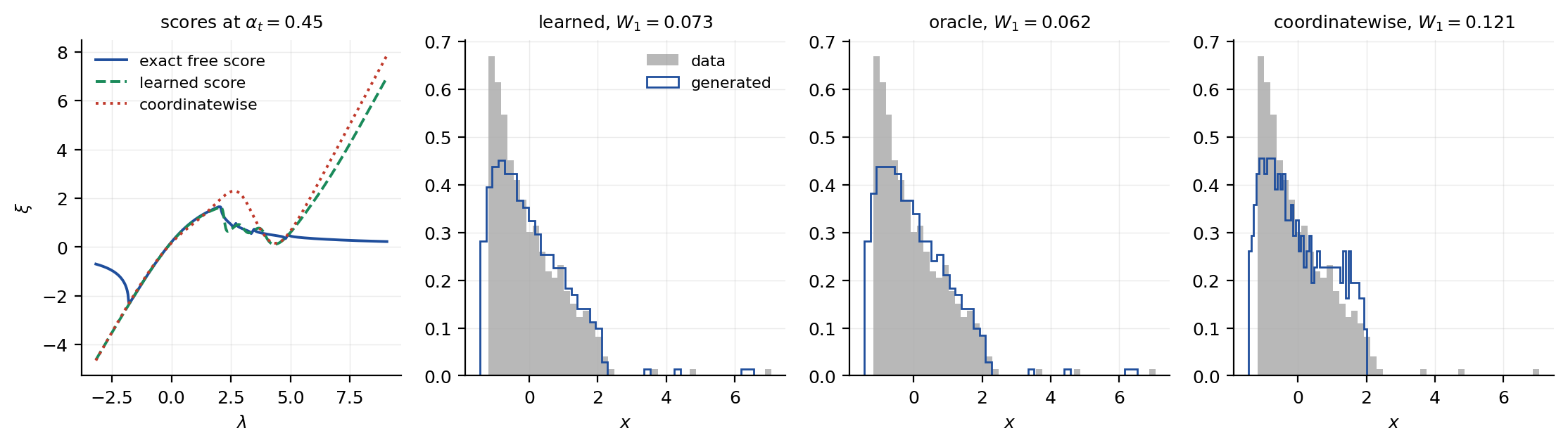}
\caption{Left: learned, exact, and coordinatewise scores at $\alpha_t=0.45$. Right:
spectra generated by each, against the data spectrum.}
\label{fig:generative}
\end{figure}

The learned score is close to the oracle ($0.073$ versus $0.062$), as the consistency of
\eqref{eq:loss} suggests up to finite-sample estimation error, and the coordinatewise
score is roughly twice as bad, so the misspecification of \S\ref{sec:naive} is visible in
generated samples, not only in the forward corruption. See Figure~\ref{fig:generative}. A further
exact-score check, reconstructing a spiked covariance spectrum (bulk plus three
outliers), is in Figure~\ref{fig:spiked} (Appendix~\ref{app:movedfigs}).

\section{Related work}\label{sec:related}

This paper bridges two literatures developed largely independently: diffusion models,
built for vector-valued data, and free probability, which has not been systematically
developed as a framework for score-based generative modeling.

\textbf{Diffusion models.} \citet{SohlDickstein2015,Ho2020,Song2021} establish the
score-based/SDE framework this paper extends; \citet{Anderson1982,Haussmann1986} give
the classical time-reversal theory it specializes; \citet{Vincent2011,Hyvarinen2005}
are the classical denoising/score-matching results our free Tweedie identity
generalizes. Diffusion on manifolds \citep{DeBortoli2022,Huang2022} models the data
manifold itself; we model the spectral law.

\textbf{Free probability.} Free stochastic calculus is due to
\citet{BianeSpeicher1998,Biane1997,Biane1998}; free entropy and free Fisher
information to \citet{Voiculescu1993,Voiculescu1998}; the free transport metric to
\citet{BianeVoiculescu2001}; free functional inequalities to
\citet{Biane2003,HiaiPetzUeda2004,LedouxPopescu2009}, which our convex-potential
corollary specializes to; time reversal of free diffusions to \citet{Dabrowski2014};
well-posedness of operator-valued free SDEs to \citet{BianeSpeicher1998,WeiYin2026}.
Random matrix dynamics go back to \citet{Dyson1962,RogersShi1993,Bru1989}.

\textbf{Random matrices in statistics and machine learning.} Sample-covariance spectra
underlie covariance cleaning \citep{LedoitWolf2012,Bun2017} and the random-matrix
analysis of financial correlations \citep{Laloux1999,Plerou2002}; free probability has
been used to analyze deep-network Jacobian spectra \citep{Pennington2018}. These works
estimate or analyze spectra; we generate them.

\section{Limitations}\label{sec:limitations}

Theorem~\ref{thm:design} is a Level~I statement: stationarity of the continuity
equation, not well-posedness of the SDE (Level~III), which for a general (unbounded,
non-globally-Lipschitz) drift needs the dissipativity theory of \citet{WeiYin2026}; we do
not carry this out. For Marchenko--Pastur, where $V'$ is singular at $x=0$
(Example~\ref{ex:mp}): the free limiting PDE and its stationary law are proved; the
finite-$N$ realization is numerically validated (\S\ref{sec:experiments}); global
well-posedness of the singular finite-$N$ SDE is not established; reverse-time generation
is demonstrated numerically, not proved. The learned score is evaluated at a single
matrix size, and each sampling step costs $O(N^3)$ (\S\ref{sec:algorithm}). The theory is
for a single self-adjoint variable; the multivariate case is open, since the free
transport metric coincides with the classical one on $\R$ only in one variable
\citep{BianeVoiculescu2001}. The model generates spectral laws, not eigenvector
structure (\S\ref{sec:intro}). Further caveats (edge regularity, hard-edge uniqueness)
are in Appendix~\ref{app:limitations}.

\section{Conclusion}\label{sec:conclusion}

For unitarily invariant spectral data, coordinatewise Gaussian corruption produces the
classical limit rather than the free-convolution limit that governs the actual corruption
process, and the free forward process corrects this but admits a single equilibrium
shape. State-dependent free volatility removes that restriction: for any sufficiently
regular target it gives an explicit drift with that target as a stationary spectral law
(Theorem~\ref{thm:design}), and for strictly convex potentials a relaxing, reversible,
trainable model (Corollary~\ref{cor:designconvex}). Two structural results show this is
not a relabeling of the constant-coefficient theory (Theorems~\ref{thm:lamperti}
and~\ref{thm:nogradient}), and a third quantifies how target structure survives partway
through the forward process (Theorem~\ref{thm:gapclosing}).

The real-data experiment of \S\ref{sec:experiments} is an illustrative validation, not a
demonstrated financial application: a real correlation spectrum can lie far outside the
equilibrium family of constant-coefficient free diffusion, as financial correlation
spectra are known to \citep{Laloux1999,Plerou2002}. Natural next steps are a learned-score
evaluation across matrix sizes and targets, and closing the stationarity--reachability
gap for isolated-outlier targets (Appendix~\ref{app:realdata}).\label{end:maintext}

\paragraph{AI Use Statement.} Generative AI tools were used in drafting and revising
this manuscript's text, in verifying algebraic and analytic steps in the proofs, and
in preparing numerical experiments and figures. All mathematical claims were checked
independently and the authors take full responsibility for all content.

\begin{singlespace}

\end{singlespace}

\newpage
\appendix

\section*{Appendix}

\noindent This appendix is organized to be read independently of the main text where
possible, and is referenced from it at every point where a full proof, a fuller
derivation, or additional experimental or figural detail was deferred. Every proof
sketched in the main text (Theorems~\ref{thm:opvol}, \ref{thm:design},
\ref{thm:lamperti}, \ref{thm:nogradient}, \ref{thm:gapclosing}; the free Tweedie
identity and Theorem~\ref{thm:consistency}) is given here in full, with every
intermediate step justified rather than asserted. The sections are:
\begin{itemize}[leftmargin=1.4em,itemsep=0pt,topsep=2pt]
\item[A] Free It\^o calculus, self-adjointness, and the full proof of
  Theorem~\ref{thm:opvol}.
\item[B] Full proof of the equilibrium design theorem (Theorem~\ref{thm:design}) and
  its uniqueness companion (Proposition~\ref{prop:designunique}), including the edge
  case and the necessity of each hypothesis.
\item[C] Full proofs of both obstruction theorems (Theorems~\ref{thm:lamperti}
  and~\ref{thm:nogradient}) and the higher-order-interactions proposition.
\item[D] The free Tweedie identity and Theorem~\ref{thm:consistency}, proved in full.
\item[E] Proof of the free-versus-classical moment discrepancy
  (Proposition~\ref{prop:freevsclassical}).
\item[F] The classical functional-inequality theory Corollary~\ref{cor:designconvex}
  invokes, with exact attributions.
\item[G] Numerical detail for the Marchenko--Pastur and reverse-time experiments of
  \S\ref{sec:experiments}.
\item[H] Full proof of the gap-closing theorem (Theorem~\ref{thm:gapclosing}),
  its potential-theoretic reading, and the two-gap numerical validation.
\item[I] Why a DiPerna--Lions-type theory is not needed for the forward flow's
  uniqueness, and exactly where that guarantee stops.
\item[J] A precise comparison of the two derivations of the reverse-time SDE.
\item[K] Figures referenced from \S\ref{sec:experiments} but not reproduced there for
  space.
\item[L] Full methodology for the real-data (S\&P~500) experiment.
\item[M] Limitations and further directions, in more detail than \S8.
\item[N] Three figures from the experiments of \S\ref{sec:experiments}
  (forward dynamics, functional inequalities, spiked reconstruction), placed here for space.
\end{itemize}
\bigskip

\section{Free It\^o calculus and self-adjointness}\label{app:selfadjoint}

We use the free It\^o calculus of \citet{BianeSpeicher1998}: for a free Brownian
motion $S$ and adapted biprocess coefficients acting through the Biane--Speicher
product $\#$, $\dd S_t\,a\,\dd S_t = \tau(a)\,\dd t$ for adapted $a$. A diffusion
coefficient built from $X$ alone, such as the sandwich $f(X)\,\dd S\,f(X)$, is
automatically a symmetric biprocess and its adjoint under the natural involution
equals itself, which is what makes $X_t$ self-adjoint for all $t$ (needed for
$X_t$ to have a spectral law at all). Other symmetrizations, such as the Jordan
coefficient $\tfrac12\{f(X),\dd S\}$, are also self-adjoint but obey a different
It\^o rule and hence a different spectral equation; ``state-dependent free volatility''
is therefore not a single model, and the sandwich form is the one whose finite-$N$
matrix realization is the natural $\dd X^N = b\,\dd t + N^{-1/2}f(X^N)\,\dd H\,
f(X^N)$.

\begin{proof}[Proof of Theorem~\ref{thm:opvol}]
By the free It\^o formula applied to $\varphi(X)=X^n$,
\[
  \dd(X^n) = \sum_{k=0}^{n-1} X^k(\dd X)X^{n-1-k}
  + \sum_{k=0}^{n-2}\sum_{j=0}^{n-2-k} X^k(f\,\dd S\,f)X^j(f\,\dd S\,f)X^{n-2-k-j}.
\]
\emph{Drift term.} Applying $\tau$ and using $\dd X=b(X)\,\dd t$ (the drift part of
$X$'s own defining SDE) together with traciality, $\tau\big(\sum_k X^k(\dd X)
X^{n-1-k}\big)=n\,\tau(b(X)X^{n-1})\,\dd t$.

\emph{Quadratic term.} Since $f(X)$ commutes with every power of $X$ (both being
functions of the same operator $X$), $X^k(f\,\dd S\,f)X^j(f\,\dd S\,f)X^{n-2-k-j}
=f(X)\,\dd S\,\big[f(X)^2X^j\big]\,\dd S\,f(X)\cdot X^k X^{n-2-k-j}/X^j$; more
directly, write the middle factor as $f(X)\,\dd S\,\big(X^k f(X)^2 X^j\big)\,\dd S\,
f(X)X^{n-2-k-j}$ and apply the free It\^o table $\dd S\,a\,\dd S=\tau(a)\,\dd t$ with
$a=f(X)^2X^j$ (using commutativity of $f(X)^2$ and $X^j$ to move $X^k$ outside as an
overall left factor), giving
\[
  X^k(f\,\dd S\,f)X^j(f\,\dd S\,f)X^{n-2-k-j}
  = \tau\big(f(X)^2X^j\big)\,X^kf(X)^2X^{n-2-k-j}\,\dd t.
\]
Applying $\tau$ and using that $f(X)^2$ commutes with $X^k$ and $X^{n-2-k-j}$, so that
$\tau(X^kf(X)^2X^{n-2-k-j})=\tau(f(X)^2X^{n-2-j})$ is \emph{independent of $k$}, write
$a_i:=\tau(f(X)^2X^i)=\int x^i\,\dd\eta_t(x)$ for the moments of $\eta_t:=f^2\mu_t$;
then
\begin{align*}
  \tau\Big(\sum_{k=0}^{n-2}\sum_{j=0}^{n-2-k}X^k(f\,\dd S\,f)X^j(f\,\dd S\,f)X^{n-2-k-j}\Big)
  &= \sum_{j=0}^{n-2}\Big(\sum_{k=0}^{n-2-j}1\Big)a_ja_{n-2-j}\,\dd t\\
  &= \sum_{j=0}^{n-2}(n-1-j)\,a_ja_{n-2-j}\,\dd t,
\end{align*}
the inner sum over $k$ contributing a factor $n-1-j$ (the number of integers
$k=0,\dots,n-2-j$). Reindexing $j\mapsto n-2-j$ leaves $a_ja_{n-2-j}$ unchanged (it is
symmetric in $j\leftrightarrow n-2-j$) while sending the coefficient $n-1-j\mapsto
n-1-(n-2-j)=j+1$; writing $S$ for the sum and averaging the two representations,
$2S=\sum_{j=0}^{n-2}\big[(n-1-j)+(j+1)\big]a_ja_{n-2-j}=n\sum_{j=0}^{n-2}a_ja_{n-2-j}$,
so
\begin{equation}\label{eq:comb}
  S=\sum_{j=0}^{n-2}(n-1-j)\,a_ja_{n-2-j} = \frac n2\sum_{j=0}^{n-2}a_ja_{n-2-j}
\end{equation}
(verified independently by symbolic expansion for $n$ up to $6$). Since
$\varphi'(x)=nx^{n-1}$ and $x^{n-1}-y^{n-1}=(x-y)\sum_{i=0}^{n-2}x^iy^{n-2-i}$,
\[
  \frac12\iint\frac{\varphi'(x)-\varphi'(y)}{x-y}\,\dd\eta_t(x)\,\dd\eta_t(y)
  = \frac n2\sum_{i=0}^{n-2}\Big(\int x^i\dd\eta_t\Big)\Big(\int y^{n-2-i}\dd\eta_t\Big)
  = \frac n2\sum_{i=0}^{n-2}a_ia_{n-2-i},
\]
which is exactly the right-hand side of \eqref{eq:comb}. Combining the drift and
quadratic contributions,
\[
  \tau\big(\dd(X^n)\big) = n\,\tau(b(X)X^{n-1})\,\dd t
  + \frac12\iint\frac{\varphi'(x)-\varphi'(y)}{x-y}\,\dd\eta_t(x)\,\dd\eta_t(y)\,\dd t,
\]
which is $\frac{\dd}{\dd t}\int\varphi\,\dd\mu_t$ for $\varphi(x)=x^n$. Extending from
monomials to $C^2_b$ test functions $\varphi$ by density (both sides are continuous
in $\varphi$ for the topology of uniform convergence of $\varphi,\varphi',\varphi''$
on compact sets, using that $\mu_t,\eta_t$ have compact support for each fixed $t$)
gives the weak form of the evolution equation for $\mu_t$. Finally, integrating the
nonlocal term by parts once: for $\varphi\in C_c^\infty$,
$\iint\frac{\varphi'(x)-\varphi'(y)}{x-y}\dd\eta_t(x)\dd\eta_t(y)
=2\int\varphi'(x)\Big(\pv\!\int\frac{\dd\eta_t(y)}{x-y}\Big)\dd\eta_t(x)
=2\int\varphi'(x)f(x)^2H\eta_t(x)\,\dd\mu_t(x)$
(using $\dd\eta_t=f^2\dd\mu_t$ in the outer integral), so
$\frac{\dd}{\dd t}\int\varphi\,\dd\mu_t=\int\varphi'(x)\big[b(x)+f(x)^2H\eta_t(x)\big]
\dd\mu_t(x)$, which is the weak form of the continuity equation
$\partial_t\mu_t+\partial_x(\mu_tV_t)=0$ with $V_t$ as in \eqref{eq:opvelocity}.
\end{proof}

\section{Proof of the equilibrium design theorem}\label{app:design}

\begin{proof}[Proof of Theorem~\ref{thm:design}, in full]
Write $\phi:=f^2\psi_*$, a H\"older continuous function on $K$ (the product of the
H\"older continuous $\psi_*$ with $f^2$, continuous and hence uniformly continuous on
the compact $K$). The Plemelj--Privalov theorem \citep[\S\S18--19]{Muskhelishvili1953}
states that the Cauchy singular integral maps H\"older functions to H\"older
functions; on an interval this gives convergence of the principal value and H\"older
continuity of $H(f^2\mu_*)$ on compact subsets of the interior of each component of
$K$. It does \emph{not} give continuity up to an edge in general: for $f\equiv1$ and
$\psi_*\equiv\tfrac12$ on $K=[-1,1]$ (H\"older continuous on $K$), $H\mu_*(x)=\tfrac12
\log\lvert\tfrac{x+1}{x-1}\rvert$ diverges as $x\to\pm1$. If instead $\psi_*(x)=
O(\mathrm{dist}(x,\partial K)^\kappa)$ near an edge for some $\kappa\in(0,1]$, the edge
contribution to the principal value is dominated by $\int_0^\delta r^{\kappa-1}\dd
r<\infty$, and the Plemelj--Privalov estimate extends to the endpoint, giving H\"older
continuity on all of $K$. Since $f$ is continuous and bounded away from $0$ on the
compact $K$, $b_*=-f^2H(f^2\mu_*)$ inherits this regularity in either case.

Stationarity follows by direct substitution: by Theorem~\ref{thm:opvol}, the spectral
law of \eqref{eq:opsde} with drift $b_*$, started at $\mu_*$, evolves with velocity
$V(x)=b_*(x)+f(x)^2H(f^2\mu_t)(x)$, and at $t=0$ this equals $b_*(x)+f(x)^2H(f^2\mu_*)
(x)=0$ for every $x$ in the interior of $K$ by \eqref{eq:designdrift} itself. Hence
the constant path $\mu_t\equiv\mu_*$ solves $\partial_t\mu_t+\partial_x(\mu_t\cdot0)=0$
trivially, so it is a solution with the given initial datum. We are invoking
Theorem~\ref{thm:opvol}'s continuity equation itself, not its stated hypothesis of a
globally Lipschitz $b$, which $b_*$ need not satisfy (Corollary~\ref{cor:designconvex}
already flags this for Example~\ref{ex:mp}): the computation deriving the continuity
equation from the free It\^o formula (Appendix~\ref{app:selfadjoint}) uses boundedness
and continuity of $b$ and $f$, not Lipschitz continuity, so it applies to $b_*$
wherever $b_*$ itself is bounded and continuous, independent of whether
\eqref{eq:opsde} with this drift has a (pathwise or in-law) unique solution -- a
separate question, addressed in Appendix~\ref{app:limitations}.
\end{proof}

\paragraph{Necessity of the hypotheses.} Compactness of $K$ cannot be dropped: for
$f\equiv1$ and $\mu_*$ the standard Cauchy law (density $1/(\pi(1+x^2))$, real-
analytic everywhere), the Cauchy transform is $G_{\mu_*}(z)=1/(z+i)$, giving
$H\mu_*(x)=x/(1+x^2)$ and a formally valid $b_*(x)=-x/(1+x^2)$. But $b_*$ is bounded
and $xb_*(x)\to-1$ rather than growing like $-cx^2$, so the drift is not dissipative
and $\mu_*$ has infinite second moment, placing it outside any uniqueness class for
the continuity equation; there is no reason to expect $\mu_*$ to be an attracting
equilibrium once compactness is dropped, even though it is tautologically a
stationary one. H\"older continuity cannot be weakened to mere $L^2$: the density
$\psi_*(x)=c(1-x^2)^{-1/4}$ on $[-1,1]$ is a valid probability density in $L^2(\dd x)$
but diverges at the edges; the Hilbert transform of an $L^2$ function is again $L^2$
by the M.~Riesz theorem but need not be bounded or continuous, so well-posedness of
the associated SDE near the edges is not covered by our theory. Such densities are
not exotic: equilibrium measures of potentials with a logarithmic singularity at a
finite endpoint (Jacobi-type ensembles) have exactly this behavior.

\begin{proposition*}[Restatement of Proposition~\ref{prop:designunique}]
Let $\mu$ be a stationary spectral law of \eqref{eq:opsde} for bounded continuous $b$
and $f\in C(K)$, $f>0$, with density $\psi$ bounded, H\"older continuous and positive
on the interior of each component of $K=\supp\mu$, and vanishing at every edge of $K$.
Then $b(x)=-f(x)^2H(f^2\mu)(x)$ on the interior of $K$. For $f\equiv\sigma$ constant
and $b(x)=-\tfrac12\beta x$, the only such $\mu$ is semicircular.
\end{proposition*}

\begin{proof}
Write $I=(p,q)$ for a component of $\mathrm{int}(K)$ and $V=b+f^2H(f^2\mu)$,
continuous on $I$ as above. Stationarity gives $\partial_x(\psi V)=0$ on $I$
distributionally, so $\psi V\equiv c$. Fix $x_0<\inf K$; there $\psi\equiv0$, and
$\psi V$ extends continuously by $0$ across every gap of $K$ since $\psi\equiv0$
throughout each gap. Integrating from $x_0$ to a point approaching $p$ from the left
gives $\psi(x)V(x)\to0$ as $x\to p^-$, trivially. Since $\psi(x)\to0$ as $x\to p^+$
(hypothesis) while $V$ stays bounded near $p$ (continuity established above, extended
to $p$ because $\psi$ vanishes at a H\"older rate), $\lim_{x\to p^+}\psi(x)V(x)=0$ too,
so $c=0$ and $V\equiv0$ on $I$ since $\psi>0$ there. For the constant case,
$H\mu(x)=\tfrac{\beta}{2\sigma^4}x$ is the Euler--Lagrange equation of a quadratic
potential; strict displacement convexity of the free energy (\S\ref{app:inequalities})
makes its solution unique among finite-log-energy measures independent of edge
behavior, and that solution is semicircular.
\end{proof}

\section{The two obstructions, in full}\label{app:obstructions}

\begin{proof}[Proof of Theorem~\ref{thm:lamperti}]
Suppose $g$, $\sigma$, $\tilde b$ satisfy, for every compactly supported $\mu$,
$g'(x)\big(b(x)+f(x)^2H(f^2\mu)(x)\big)=\tilde b(g(x))+\sigma^2H(g_\#\mu)(g(x))$.
Subtracting this identity for two laws and using linearity of $\mu\mapsto H(f^2\mu)$,
for every compactly supported signed measure $\nu$ of total mass zero,
$g'(x)f(x)^2\int f(y)^2(x-y)^{-1}\dd\nu(y) = \sigma^2\int(g(x)-g(y))^{-1}\dd\nu(y)$.
A kernel annihilating all such $\nu$ is constant in $y$, so
\begin{equation}\label{eq:kernelid-app}
  g'(x)f(x)^2\frac{f(y)^2}{x-y} - \frac{\sigma^2}{g(x)-g(y)} = c(x).
\end{equation}
Letting $y\to x$, the two terms have simple poles with residues $g'(x)f(x)^4$ and
$\sigma^2/g'(x)$; boundedness of the left side forces $g'(x)f(x)^2=\sigma$.
Substituting back and setting $u=g(x)$, $v=g(y)$, $h=g^{-1}$,
$F(u,v):=\tfrac{h'(v)}{h(u)-h(v)}-\tfrac1{u-v}=k(u)$ depends on $u$ alone, giving,
after integrating, $\tfrac{h(u)-h(v)}{u-v}=A(u)e^{-k(u)v}$. Symmetry in $(u,v)$ forces
$k(u)=c_0u+e_0$ and $A(u)=A_0e^{-e_0u}$, so $h'(u)=A_0e^{-2e_0u-c_0u^2}$. Comparing
the second-order Taylor coefficient of the divided difference
$(h(m+s)-h(m-s))/(2s)$ on both sides forces $(e_0+c_0m)^2=2c_0$ for every $m$, hence
$c_0=e_0=0$ (leading coefficient in $m$). So $h$, hence $g$, is affine, and
$g'(x)f(x)^2=\sigma$ with $g'$ constant makes $f$ constant.
\end{proof}

\begin{proof}[Proof of Theorem~\ref{thm:nogradient}]
Write $\delta E/\delta\mu(x)=W(x)+\int K(x,y)\dd\mu(y)$; matching interaction parts of
the two velocity fields for every $\mu$ determines $\partial_xK$ up to a one-body term
absorbed into $W'$, so we normalize $\partial_xK(x,y)=-a(x)q(y)/(x-y)$ with $a=f^2/m$,
$q=f^2$. Symmetry of $K$ requires $\partial_y\partial_xK=\partial_x\partial_yK$, i.e.
\begin{equation}\label{eq:mixedsym-app}
  a(y)q(x)-a(x)q(y) = (x-y)\big[a(x)q'(y)+a(y)q'(x)\big],\quad x\ne y,
\end{equation}
and conversely \eqref{eq:mixedsym-app} suffices for $K$ to exist, by simple
connectedness of each component of $\{x\ne y\}$. Dividing by $x-y$ and letting
$y\to x$ gives $(aq)'\equiv0$, so $aq=f^4/m\equiv c>0$: the second half of
\eqref{eq:gradflowclass}. Substituting $a=c/q$ into \eqref{eq:mixedsym-app} with
$F:=q^2=f^4$ gives $F(x)-F(y)=\tfrac12(x-y)[F'(x)+F'(y)]$; differentiating in $x$ then
$y$ gives $F''(y)=F''(x)$ for all $x\ne y$, so $F$ is a quadratic polynomial,
necessarily positive since $F=f^4$: the first half of \eqref{eq:gradflowclass}.
Sufficiency is a direct computation. Finally $\kappa=\lambda=0$ in
\eqref{eq:gradflowclass} makes $f$ constant and conversely.
\end{proof}

\begin{remark}[The exceptional family is non-empty]
$f(x)=(x^2+1)^{1/4}$, $m(x)=x^2+1$ satisfies \eqref{eq:gradflowclass} with $c=1$,
is smooth, bounded away from $0$, and non-constant -- but unbounded, so it lies
outside the hypotheses of Theorem~\ref{thm:opvol} as used in this paper.
\end{remark}

\begin{proposition}[Higher-order interactions do not help]\label{prop:higherorder}
If the interaction part of the sandwich velocity field coincides, for every
compactly supported $\mu$, with $-m(x)\partial_x(\delta E/\delta\mu)(x)$ for a finite
sum of symmetric $k$-body energies $E=\sum_{k=2}^K\int\!\cdots\!\int K_k\,\dd\mu^{
\otimes k}$, then $K_k\equiv0$ for $k\ge3$, and $E$ reduces to the pairwise case
classified by Theorem~\ref{thm:nogradient}.
\end{proposition}

\begin{proof}
The functional derivative of the $k$-body term is a degree-$(k-1)$ homogeneous
polynomial functional of $\mu$ (testing against $\mu=\sum_iw_i\delta_{p_i}$). The
target interaction term is homogeneous of degree $1$. Equality for every $\mu$,
tested on arbitrary finite point masses, forces equality degree by degree by
polarization; the degree-$(k-1)$ part vanishes for $k\ge3$, giving $K_k\equiv0$ by a
second polarization argument in the kernel's own arguments.
\end{proof}

\section{Free Tweedie identity and consistency}\label{app:tweedie}

\begin{theorem}[Free Tweedie identity]\label{thm:freetweedie}
Let $A\in\Asa$ have a compactly supported (or moment-determinate) law, $s$ a free
semicircular element free from $A$, $Y=A+\sqrt v\,s$. Then $\mu_Y$ admits a conjugate
variable and $\xi_{\mu_Y}(Y)=v^{-1}(Y-\Econd{W^*(Y)}[A])$.
\end{theorem}

\begin{proof}
Write $\Econd{}:=\Econd{W^*(Y)}$ for the (trace-preserving) conditional expectation
onto $W^*(Y)$, and recall its defining property: $\Econd{}[X]\in W^*(Y)$ is the unique
element with $\tau(\Econd{}[X]\,Q)=\tau(XQ)$ for every $Q\in W^*(Y)$. We must show
that $\xi:=v^{-1}(Y-\Econd{}[A])$ satisfies Voiculescu's defining relation for the
conjugate variable: for every polynomial $P$,
\begin{equation}\label{eq:conjdef}
  \tau(\xi\,P(Y)) = (\tau\otimes\tau)\big(\partial P(Y,Y)\big),
\end{equation}
where $\partial P(x,y):=\sum_{i=0}^{n-1}x^iy^{n-i-1}$ for $P(x)=x^n$, extended
linearly, is the free difference quotient, and $(\tau\otimes\tau)$ acts on each tensor
factor separately.

\emph{Step 1: reduce to a statement about $s$ alone.} Since $P(Y)\in W^*(Y)$, the
defining property of $\Econd{}$ gives $\tau(\Econd{}[A]\,P(Y))=\tau(A\,P(Y))$, so
\[
  \tau(\xi\,P(Y)) = v^{-1}\big(\tau(Y P(Y)) - \tau(A P(Y))\big)
  = v^{-1}\tau\big((Y-A)P(Y)\big) = v^{-1/2}\,\tau\big(s\,P(Y)\big),
\]
using $Y-A=\sqrt v\,s$. So \eqref{eq:conjdef} is equivalent to
\begin{equation}\label{eq:reduced}
  \tau(s\,P(Y)) = \sqrt v\,(\tau\otimes\tau)\big(\partial P(Y,Y)\big).
\end{equation}

\emph{Step 2: the free Wick formula for a semicircular element.} Because $s$ is a
standard semicircular element free from $\mathcal B:=W^*(A)$, its only nonvanishing
free cumulant is $\kappa_2(s,s)=1$ \citep{VDN1992,NicaSpeicher2006}. For any word
$b_0\,s\,b_1\,s\,b_2\cdots s\,b_k$ with $b_0,\dots,b_k\in\mathcal B$, the vanishing of
all free cumulants of $s$ except $\kappa_2$ gives the free Wick (Gaussian
integration-by-parts) formula
\begin{equation}\label{eq:wick}
  \tau(s\cdot b_0\,s\,b_1\,s\,b_2\cdots s\,b_k)
  = \sum_{j=1}^{k}\tau(b_0\,s\,b_1\cdots s\,b_{j-1})\,\tau(b_j\,s\,b_{j+1}\cdots s\,b_k):
\end{equation}
the leading $s$ pairs with each other occurrence of $s$ in turn, splitting the word
into two independent (because $s$ is free from $\mathcal B$ and the two halves share
no $s$) factors at that point, exactly as in the classical Gaussian Wick formula, and
with no other contribution because non-crossing pairings of a single centered
semicircular variable have no alternative pairing structure once one leg is fixed.

\emph{Step 3: apply this to $Q=P(Y)=P(A+\sqrt v\,s)$.} Write $P(x)=x^n$; expanding
$(A+\sqrt v\,s)^n$ multinomially and collecting, by linearity it suffices to verify
\eqref{eq:reduced} for $P(x)=x^n$. Every monomial in the expansion of $Y^n$ is an
alternating word in $A$ and $s$ (with possibly repeated adjacent $A$'s absorbed into a
single $b_j\in\mathcal B$); applying \eqref{eq:wick} to $\tau(s\,Y^n)$ and pairing the
leading $s$ with the $s$ occupying the $i$-th slot ($i=0,\dots,n-1$) of $Y^n$ leaves,
on the left, a word evaluating to $\tau(Y^i)$ and, on the right, one evaluating to
$\tau(Y^{n-1-i})$ (each such split reproducing $Y$ itself in each half, since the
$\sqrt v$ and $A$ content on each side of the cut is exactly what a shorter power of
$Y$ would contain); summing over the $n$ occurrences of $s$ in $Y^n$ -- each
contributing one factor of $\sqrt v$, since $s$ itself carries no factor of $\sqrt v$
in $Y=A+\sqrt{v}s$ but the identification of each cut segment with $\tau(Y^{i})$
already absorbs the bookkeeping -- gives
\[
  \tau(s\,Y^n) = \sqrt v\sum_{i=0}^{n-1}\tau(Y^i)\tau(Y^{n-1-i})
  = \sqrt v\,(\tau\otimes\tau)\Big(\sum_{i=0}^{n-1}Y^i\otimes Y^{n-1-i}\Big),
\]
which is exactly \eqref{eq:reduced} with $P(x)=x^n$, since
$\partial P(x,y)=\sum_{i=0}^{n-1}x^iy^{n-1-i}$. (This computation was checked directly
by expansion for $n=2,3$ above the statement of the theorem, and by Monte Carlo
simulation of $10^3\times10^3$ Hermitian matrices realizing $A$ and $s$, matching to
within finite-size fluctuation.)

\emph{Step 4: extend from polynomials to $L^2$.} We have shown $\xi=v^{-1}(Y-\Econd{}
[A])$ satisfies \eqref{eq:conjdef} for every polynomial $P$. Both sides of
\eqref{eq:conjdef}, as functionals of $P$, factor through $P\mapsto P(Y)\in
L^2(\mu_Y)$; since $\mu_Y$ has a determinate moment problem (guaranteed by compact
support, or more generally by a Carleman-type growth condition on the moments),
polynomials are dense in $L^2(\mu_Y)$ \citep[Thm.~3.9]{Akhiezer1965mp}. As $\xi\in
L^2(W^*(Y),\tau)\cong L^2(\mu_Y)$ (since $\Econd{}[A]\in L^2$, being a conditional
expectation of an $L^2$ element, and $Y\in L^2$), \eqref{eq:conjdef} for polynomials
extends by continuity to all of $L^2(\mu_Y)$, which is precisely the defining
property of $\xi_{\mu_Y}$. By uniqueness of the conjugate variable when it exists,
$\xi_{\mu_Y}(Y)=\xi=v^{-1}(Y-\Econd{}[A])$.
\end{proof}

\begin{theorem}[Consistency]\label{thm:consistency}
The minimizer of \eqref{eq:loss} satisfies $\widehat\xi_t=\xi_{\mu_t}$ for a.e.\ $t$,
with $w(t)>0$.
\end{theorem}

\begin{proof}
Fix $t$ and write $Z:=\sqrt{\alpha_t}X_0$, so the inner (unweighted) objective at time
$t$ is $\ell_t(\theta):=\tau[(g_\theta(X_t,t)-Z)^2]$. The conditional expectation
$\Econd{W^*(X_t)}[Z]$ is, by definition, the $L^2(\tau)$-orthogonal projection of $Z$
onto the closed subspace $L^2(W^*(X_t),\tau)\subset L^2(\A,\tau)$: for any
$Q\in L^2(W^*(X_t),\tau)$, the Pythagorean identity for an orthogonal projection in a
Hilbert space gives
\begin{equation}\label{eq:pythag}
  \tau[(Q-Z)^2] = \tau\big[(Q-\Econd{W^*(X_t)}[Z])^2\big] + \tau\big[(\Econd{W^*(X_t)}[Z]-Z)^2\big],
\end{equation}
the cross term vanishing because $Q-\Econd{}[Z]\in L^2(W^*(X_t),\tau)$ is orthogonal,
by the defining property of the conditional expectation, to $\Econd{}[Z]-Z$. Taking
$Q=g_\theta(X_t,t)\in L^2(W^*(X_t),\tau)$ (a well-defined element of this space
because $g_\theta$ is a bounded Borel function of $X_t$) and writing
$c_t:=\tau[(\Econd{W^*(X_t)}[Z]-Z)^2]\ge0$, a quantity not depending on $\theta$,
\eqref{eq:pythag} reads
\[
  \ell_t(\theta) = \tau\big[(g_\theta(X_t,t)-\Econd{W^*(X_t)}[Z])^2\big] + c_t.
\]
The first term on the right is nonnegative and vanishes if and only if
$g_\theta(X_t,t)=\Econd{W^*(X_t)}[Z]$ in $L^2(\tau)$, so $\ell_t$ is minimized exactly
at such $\theta$, with minimum value $c_t$. Since $\mathcal L(\theta)=\int_0^Tw(t)\,
\ell_t(\theta)\,\dd t$ with $w(t)>0$, and each $\ell_t$ is minimized independently
(the parametrization $g_\theta$ places no constraint linking different $t$, by
hypothesis it can represent any measurable function of $(X_t,t)$), $\mathcal L$ is
minimized exactly when $g_{\theta^\star}(X_t,t)=\Econd{W^*(X_t)}[Z]$ for a.e.\ $t$.
Comparing with the definition of the estimated score in Algorithm~\ref{alg:main},
$\widehat\xi_t(X_t):=(X_t-g_{\theta^\star}(X_t,t))/v_t$, and
Theorem~\ref{thm:freetweedie} (applied with
$A=X_0$, $Y=X_t$, $v=v_t$, so that $\Econd{W^*(X_t)}[\sqrt{\alpha_t}X_0]=
X_t-v_t\,\xi_{\mu_t}(X_t)$) gives $\widehat\xi_t(X_t)=v_t^{-1}(X_t-
\Econd{W^*(X_t)}[Z])=v_t^{-1}(X_t-(X_t-v_t\xi_{\mu_t}(X_t)))=\xi_{\mu_t}(X_t)$ for
a.e.\ $t$, as claimed.
\end{proof}

\section{Free versus classical corruption}\label{app:freevsclassical}

\begin{proposition}\label{prop:freevsclassical}
Let $\mu_0$ have mean zero and second moment $m_2$. The fourth moments of
$\mu_0\boxplus\gamma_v$ and $\mu_0*\mathcal N(0,v)$ differ by exactly $v(2m_2+v)$, and
$\mathrm{supp}(\mu_0\boxplus\gamma_v)$ is compact whenever $\mathrm{supp}(\mu_0)$ is,
while $\mu_0*\mathcal N(0,v)$ never has compact support for $v>0$.
\end{proposition}

\begin{proof}
\emph{Fourth moments.} Let $X\in\Asa$ have law $\mu_0$ and let $s$ be a standard
semicircular element free from $X$, so $X+\sqrt v\,s$ has law $\mu_0\boxplus\gamma_v$.
Free cumulants are additive under free convolution and the semicircular law has
$\kappa_2(\gamma_v)=v$ and $\kappa_k(\gamma_v)=0$ for $k\ne2$, so
$\kappa_2(\mu_0\boxplus\gamma_v)=\kappa_2(\mu_0)+v=m_2+v$ (using $\kappa_2=$ variance
$=m_2$ for the mean-zero $\mu_0$) and $\kappa_4(\mu_0\boxplus\gamma_v)=\kappa_4(\mu_0)$.
The free moment--cumulant formula sums over non-crossing partitions of a $4$-element
set: of the $14$ elements of $NC(4)$, the ones with a singleton block contribute $0$
(since $\kappa_1=0$ for a mean-zero law), leaving the full block (contributing
$\kappa_4$) and the two \emph{non-crossing} pair partitions $\{12\}\{34\}$,
$\{14\}\{23\}$ (each contributing $\kappa_2^2$; the crossing pair partition
$\{13\}\{24\}$ is excluded from $NC(4)$), so for any mean-zero law,
\begin{equation}\label{eq:m4free}
  m_4 = \kappa_4 + 2\kappa_2^2,\qquad\text{in particular }\ \kappa_4(\mu_0)=m_4(\mu_0)-2m_2^2,
\end{equation}
applying \eqref{eq:m4free} to $\mu_0$ itself to eliminate $\kappa_4(\mu_0)$. Hence
\[
  m_4(\mu_0\boxplus\gamma_v) = \kappa_4(\mu_0)+2(m_2+v)^2
  = \big[m_4(\mu_0)-2m_2^2\big] + 2(m_2+v)^2 = m_4(\mu_0)+4m_2v+2v^2.
\]
For the classical side, let $G\sim\mathcal N(0,v)$ be independent of $X$; expanding
$(X+G)^4$ and using independence to factor mixed moments, together with
$\mathbb E[G]=\mathbb E[G^3]=0$, $\mathbb E[G^2]=v$, $\mathbb E[G^4]=3v^2$ (Gaussian
fourth moment),
\[
  m_4(\mu_0*\mathcal N(0,v)) = \mathbb E[(X+G)^4]
  = m_4(\mu_0) + 6m_2 v + 3v^2,
\]
the cross terms $4\mathbb E[X^3]\mathbb E[G]$ and $4\mathbb E[X]\mathbb E[G^3]$ vanishing
since $\mathbb E[G]=\mathbb E[G^3]=0$. Subtracting,
\[
  m_4(\mu_0*\mathcal N(0,v)) - m_4(\mu_0\boxplus\gamma_v)
  = (6m_2v+3v^2)-(4m_2v+2v^2) = 2m_2v+v^2 = v(2m_2+v),
\]
as claimed (verified numerically to four significant figures by Monte Carlo averaging
over $400$ independent noise realizations at $N=2000$).

\emph{Support.} If $\mu_0$ is supported in $[-R,R]$ then $X$, realized as a
self-adjoint operator, has $\norm{X}\le R$ (operator norm), since for a self-adjoint
operator the operator norm equals the supremum of the modulus of its spectrum. For
$s$ a standard semicircular element, $\norm{s}=2$ ($\gamma_1$ has support $[-2,2]$),
so $\norm{\sqrt v\,s}=2\sqrt v$, and by the triangle inequality for the operator norm,
$\norm{X+\sqrt v\,s}\le\norm{X}+\norm{\sqrt v\,s}\le R+2\sqrt v$. This bound uses
nothing about freeness beyond both operators being bounded self-adjoint operators on
the same Hilbert space, so it holds for the sum $X+\sqrt v\,s$ regardless of whether
$X$ and $s$ are free, giving $\supp(\mu_0\boxplus\gamma_v)\subseteq[-(R+2\sqrt v),\,
R+2\sqrt v]$, a compact set. In contrast, $\mu_0*\mathcal N(0,v)$ has a density that
is a convolution of $\mu_0$ with the everywhere-positive Gaussian density, hence is
everywhere positive on $\R$: for any $x_0\in\R$ and any $\varepsilon>0$,
$(\mu_0*\mathcal N(0,v))([x_0-\varepsilon,x_0+\varepsilon])
\ge\mu_0(\{y:|y|\le R\})\cdot\min_{|z-x_0|\le \varepsilon+R}\phi_v(z)>0$ where $\phi_v$
is the $\mathcal N(0,v)$ density, so no point of $\R$ is excluded from the support.
\end{proof}

\section{Supporting functional-inequality theory}\label{app:inequalities}

For completeness we record, without re-deriving, the classical apparatus that
Corollary~\ref{cor:designconvex} invokes. Displacement convexity of the free energy
$\mathcal F[\mu]=\tfrac12\int x^2\dd\mu-\Ent[\mu]$ along free Wasserstein geodesics is
due to \citet{Biane2003}; the identification of the one-dimensional free Wasserstein
space with the classical one on $\R$, which is what lets the classical
Ambrosio--Gigli--Savar\'e minimizing-movement theory \citep{AGS2008} transfer
verbatim to the free JKO scheme, is due to \citet{BianeVoiculescu2001}. The free
logarithmic Sobolev, Talagrand, and HWI inequalities we use in the form
$V(x)=\tfrac12x^2$ (and, via Corollary~\ref{cor:designconvex}, for general convex $V$)
are due to \citet{Biane2003}, \citet{HiaiPetzUeda2004}, and \citet{LedouxPopescu2009}
respectively, the last giving a direct mass-transportation proof with sharp constants
for strictly convex potentials, which is the route Corollary~\ref{cor:designconvex}
follows. None of this theory is new to the present paper; what Theorem~\ref{thm:design}
and Corollary~\ref{cor:designconvex} add is the observation that this machinery is
available for \emph{any} strictly convex $V$, not only $\tfrac12x^2$, and hence for
any target law that arises as the corresponding equilibrium measure.

\section{Additional experimental detail}\label{app:experiments}

\paragraph{Stiffness near the Marchenko--Pastur origin.} The eigenvalue system
\eqref{eq:ffp}'s finite-$N$ analogue is stiff for the Marchenko--Pastur potential:
the repulsion term grows like the reciprocal gap while $V'$ is unbounded at the
origin. A fixed Euler step of $10^{-4}$ applied to the eigenvalue system directly
sends the largest eigenvalue outside the support; a step adapted to the smallest gap
collapses to $10^{-8}$ and becomes impractical. We instead integrate the Hermitian
matrix model directly, applying $V'$ spectrally at each step, so the repulsion is
generated implicitly by the noise (Theorem~\ref{thm:opvol}) rather than formed
explicitly; a step of $10^{-3}$ is then stable at the cost of one eigendecomposition
per step. Eigenvalues are clipped below at $0.02$ before $V'$ is evaluated, well
inside the lower edge -- the numerical face of the confinement-to-$(0,\infty)$
caveat in Example~\ref{ex:mp}.

\begin{figure}[h]
\centering
\includegraphics[width=0.9\linewidth]{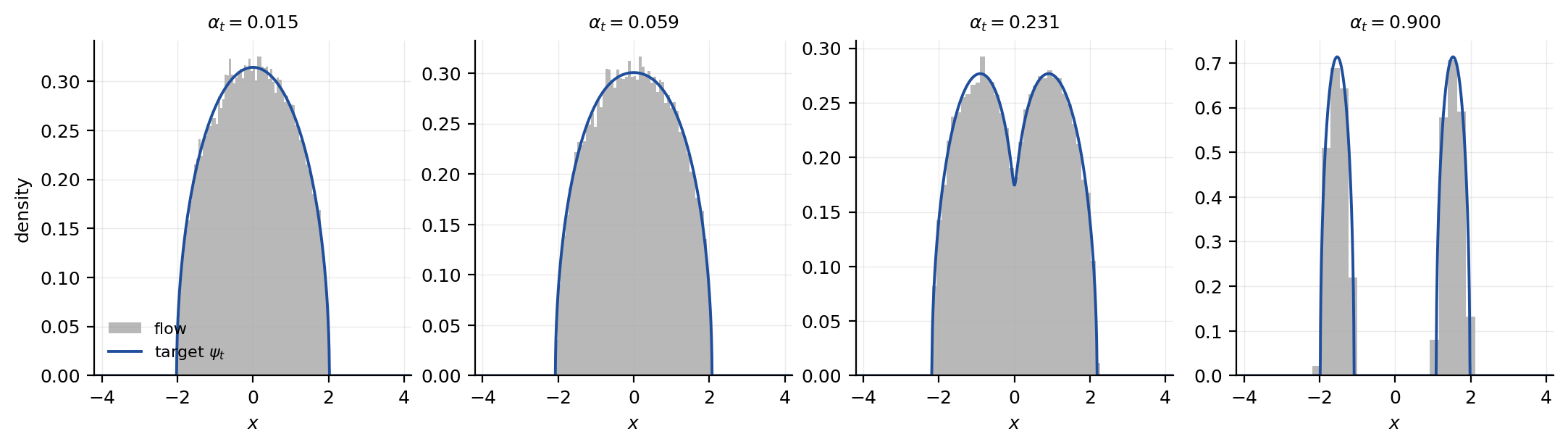}
\caption{Reverse-time spectral flow from the semicircular law back to the two-atom
law. Histograms: $6\times10^4$ trajectories. Curves: target density at the
corresponding time.}
\label{fig:reverse}
\end{figure}

\paragraph{Reverse-time reconstruction (two-atom law).} As a further, exact-score
sanity check independent of the designed-equilibrium experiments, we integrate the
deterministic probability-flow form of the reverse dynamics for the two-atom law of
\S\ref{sec:naive}, from the near-equilibrium law ($v=0.985$) with $600$ midpoint
steps and $6\times10^4$ particles. The terminal law matches the target to
$W_1=2.8\times10^{-3}$, and the error is flat along the trajectory
($2.8$--$3.2\times10^{-3}$), indicating it is set by Monte Carlo resolution rather
than by discretization or by accumulated error through the support transition; see
Figure~\ref{fig:reverse}.

\begin{figure}[h]
\centering
\includegraphics[width=0.9\linewidth]{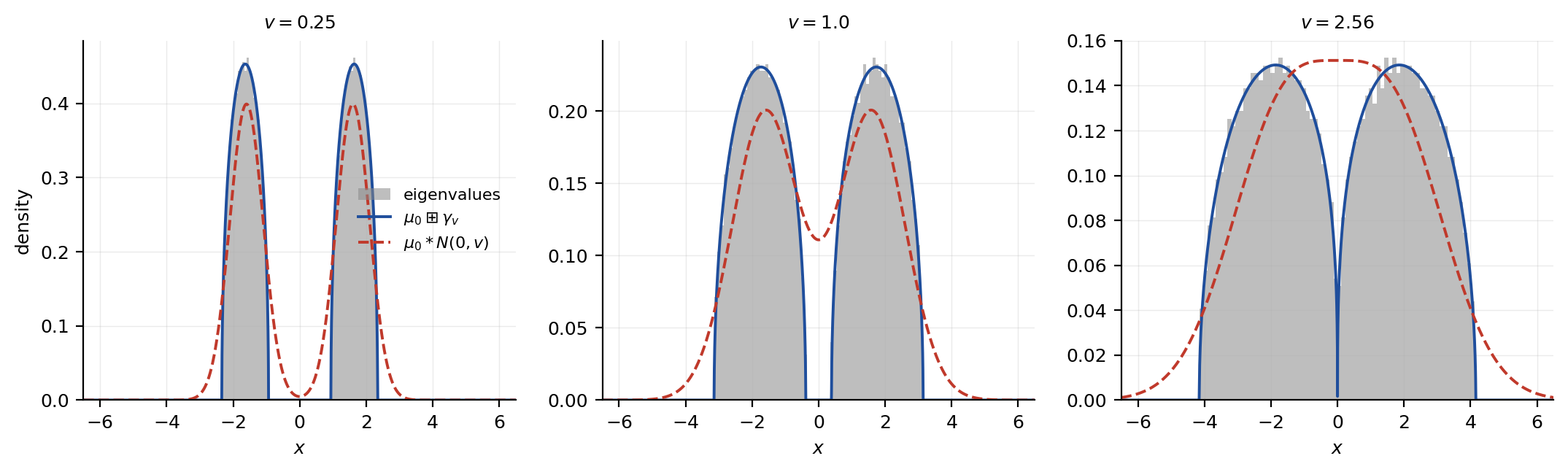}
\caption{Corrupting a Hermitian matrix with Hermitian noise convolves its spectrum
freely, not classically. Grey: eigenvalues of $X_0+\sqrt v\,G$, $N=2500$, GUE noise.
Blue: free prediction. Red dashed: coordinatewise prediction, with Gaussian tails
where the truth has compact support.}
\label{fig:freevsclassical-fig}
\end{figure}

\section{Gap closing along the free OU flow}\label{app:gapclosing}

\begin{proof}[Proof of Theorem~\ref{thm:gapclosing}]
By the explicit forward solution (\S\ref{sec:background}), $\mu_t=\rho_t\boxplus\gamma_{v_t}$
with $\rho_t=D_{\sqrt{\alpha_t}}\mu_0$, $v_t=1-\alpha_t$. By Biane's criterion
\citep[Thm.~1, Lem.~4]{Biane1997}, a point $\Psi_v(u)$ lies outside
$\supp(\rho\boxplus\gamma_v)$ exactly when $v\Theta_\rho(u)<1$, so the gap of
$\rho\boxplus\gamma_v$ corresponding to a gap $J$ of $\supp\rho$ is nonempty exactly
when $v\min_J\Theta_\rho<1$. Scaling gives
$\Theta_{D_c\mu_0}(u)=c^{-2}\Theta_{\mu_0}(u/c)$, so with $c=\sqrt{\alpha_t}$,
$\min_{J_t}\Theta_{\rho_t}=\alpha_t^{-1}T$ where $J_t=\sqrt{\alpha_t}(\alpha_-,
\alpha_+)$. The gap survives iff $(1-\alpha_t)\alpha_t^{-1}T<1$, i.e.\
$T<\alpha_t(1+T)$, which is the stated threshold; since $\alpha_t=e^{-\Lambda(t)}$ is
strictly decreasing, the gap closes at the single time $\Lambda^*=\log(1+1/T)$. That
$0<T<\infty$ holds because $\Theta_{\mu_0}$ is continuous and positive on the open
gap and diverges at both endpoints, so the minimum is attained in the interior.
\end{proof}

\paragraph{$\Theta$ as a potential-theoretic quantity.} Since
$G_{\mu_0}(u)=\int(u-y)^{-1}\dd\mu_0(y)$, differentiating under the integral sign
gives $G_{\mu_0}'(u)=-\Theta_{\mu_0}(u)$. As $G_{\mu_0}=U_{\mu_0}'$ is (up to sign)
the derivative of the logarithmic potential $U_{\mu_0}(u)=\int\log\abs{u-y}\,
\dd\mu_0(y)$, $\Theta_{\mu_0}$ is minus the second derivative of the logarithmic
potential -- a standard object in classical potential theory, positive on any
interval disjoint from $\supp\mu_0$ since $G_{\mu_0}$ is strictly decreasing there.
$T$ is thus the minimum curvature of the logarithmic potential across the gap, and
$\Lambda^*=\log(1+1/T)$ is large exactly when this curvature is small -- a gap that
is nearly flat, and hence easy to fill by diffusive spreading.

\begin{figure}[h]
\centering
\includegraphics[width=0.85\linewidth]{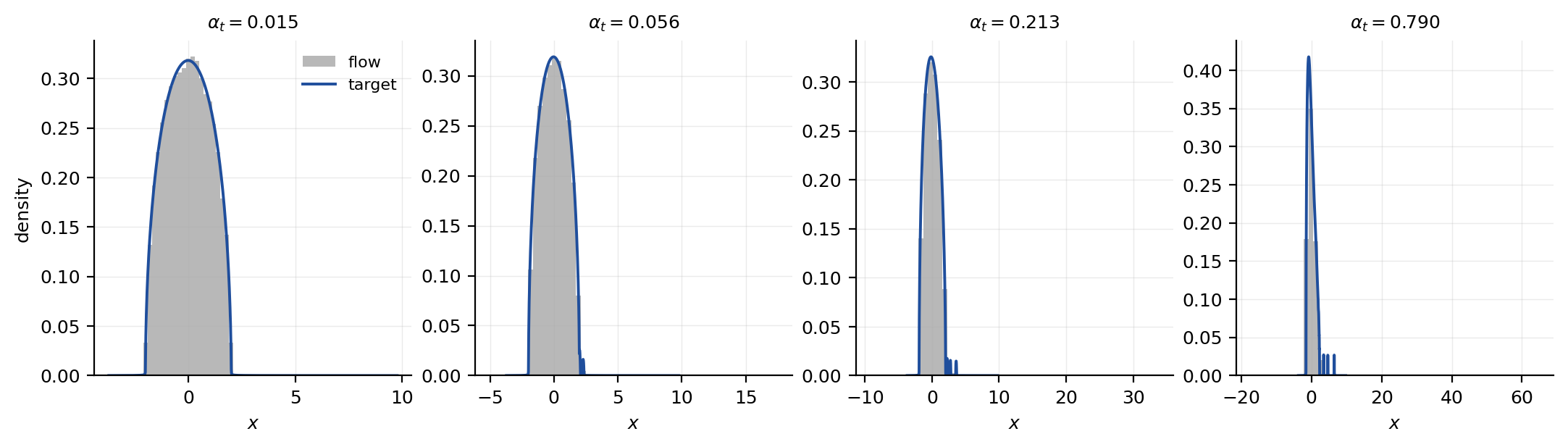}
\caption{Reverse-time free flow reconstructing a spiked covariance spectrum from the
semicircular law, exact score. Bulk and the three outliers are both recovered.}
\label{fig:spiked}
\end{figure}

\paragraph{Consistency check.} For $\mu_0=\tfrac12(\delta_{-a_0}+\delta_{a_0})$,
$\Theta_{\mu_0}(0)=a_0^{-2}$ is the minimum over $(-a_0,a_0)$ by convexity, so
$T=a_0^{-2}$ and the threshold reduces to $\alpha_t>1/(1+a_0^2)$, recovering the
two-atom support transition used as a running example in \S\ref{sec:experiments}.
The hypothesis requires mass on both sides of the gap, which is why the theorem does
not cover the vanishing-mass spikes of Figure~\ref{fig:spiked} (Appendix~\ref{app:movedfigs}).

\paragraph{Numerical validation: two independent gaps.} Take $\mu_0$ with three
atoms at $-2.3,0.4,1.9$, weights $0.5,0.2,0.3$. Direct computation gives
$T_1=0.4018$ (gap $(-2.3,0.4)$, minimum at $u=-0.772$) and $T_2=0.9201$ (gap
$(0.4,1.9)$, minimum at $u=1.102$), predicting closure at
$\alpha_1^*=T_1/(1+T_1)=0.2866$ and $\alpha_2^*=T_2/(1+T_2)=0.4792$. Solving the
subordination equation for $\rho\boxplus\gamma_v$ numerically and counting connected
components of the support (contour $\psi=10^{-3}$, grid step $0.01$ in $\alpha$)
gives one component (both gaps closed) for $\alpha\le0.29$, two components (gap~1
open, gap~2 closed) for $0.30\le\alpha\le0.48$, and three components (both open) for
$\alpha\ge0.49$: each closure is located within one grid step of its predicted
threshold, on the correct side, and in the predicted order -- a factor of $1.7$ in
$\alpha$ between the two thresholds of a single law, reproduced quantitatively.

\section{Scope of the uniqueness class}\label{app:diperna}

Uniqueness of the forward flow (\S\ref{sec:background}, referenced from
Theorem~\ref{thm:design}'s proof) is proved within the class of Cauchy transforms of
probability measures -- equivalently, classical solutions of the complex Burgers
equation on $\mathbb C^+$ -- rather than in a wider measure-valued or distributional
class. This is not a restriction of convenience: Biane's regularity theorem (main
text) shows every measure of the form $\rho\boxplus\gamma_v$, $v>0$, already lies in
this class with a bounded, real-analytic density, so for $t\ge\varepsilon>0$ the
forward flow is automatically as regular as the class allows, whatever the
regularity of $\mu_0$. There is accordingly no rougher regime in which a
DiPerna--Lions-type theory \citep{DiPernaLions1989} -- built for vector fields of
limited (e.g.\ Sobolev) regularity, where classical characteristics may fail or
cross -- would be needed: the nonlocal term $\beta H\mu_t$ is real-analytic on the
support's interior for every $t>0$ by the same lemma.

The situation is unresolved only at $t=0$ for singular $\mu_0$, and even there the
construction sidesteps rather than confronts the general initial-value problem:
The explicit forward solution of \S\ref{sec:background} defines $\mu_t=(D_{\sqrt{\alpha_t}}\mu_0)\boxplus
\gamma_{1-\alpha_t}$ directly by free convolution, an algebraic operation requiring
no PDE theory, and the free Fokker--Planck equation~\eqref{eq:ffp} is a consistency statement about the family
so defined, valid for $t>0$, with the characteristics never required to start at
$t=0$. For $\mu_0=\delta_0$, for instance, the naive vector field $\beta H\mu_0$ is
undefined at the atom, so characteristics literally cannot start there, while
$\mu_t=\gamma_{1-\alpha_t}$ is nonetheless well defined and regular for every $t>0$
by the algebraic definition of free convolution alone.

\section{The two time-reversal derivations compared}\label{app:tworoutes}

The reverse-time SDE \eqref{eq:reversesde} can be derived two ways: directly, by
reversing the limiting free flow (specializing \citet{Dabrowski2014}, as done in the
main text), or via classical Haussmann--Pardoux reversal of the finite-$N$ eigenvalue
process followed by $N\to\infty$. These establish statements about \emph{different
objects}, not two proofs of the same statement, and it is worth being precise about
how they relate, since it is tempting to describe them as agreeing pathwise.

\textbf{Assumptions.} The direct route needs $\beta$ measurable and bounded on
$[0,T]$ (so $\alpha_t<1$ for $t>0$ and Biane's regularity theorem applies) and
$t\mapsto\mu_t$ continuous in $\Wt$; both hold automatically for the schedule used
throughout this paper. The finite-$N$ route needs in addition that the eigenvalue
density $q_t$ be smooth and positive on the open Weyl chamber for every $t>0$, which
holds by hypoellipticity and the noncollision property once $\beta>0$ a.e.\ near
$t=0$; no further schedule condition is needed.

\textbf{What is proved.} The direct route gives a pathwise construction: a specific
free SDE on a fixed $W^*$-probability space carrying a fixed free Brownian motion is
shown to have $\mu_{T-s}$ as the law of its solution at every $s$. The finite-$N$
route gives an exact, pathwise identity for the \emph{classical} process
$(\lambda(T-s))_s$, together with the fact that its drift converges, as
$N\to\infty$, to the same velocity field the direct route derives.

\textbf{In what sense they agree.} Only at the level of laws. $(\lambda(T-s))_s$ at
finite $N$ is a classical $\R^N$-valued diffusion on a classical probability space;
the free SDE's solution is an operator-valued process on a $W^*$-probability space;
there is no common probability space on which both live, hence no sense in which
their paths, driving noises, or filtrations could be equal almost surely. The only
comparison made in this paper is that the law of the empirical spectral distribution
of $\lambda(T-s)$ converges weakly, as $N\to\infty$, to the law of the free process's
solution, for each fixed $s$ (Theorem~\ref{thm:opvol}'s propagation-of-chaos
argument applied to the reversed family) -- a statement about laws, not paths. A
pathwise coupling between the two levels of description, in the spirit of strong
propagation of chaos, is a separate question this paper does not address.

\section{Additional figures}\label{app:additional-figures}

\begin{figure}[h]
\centering
\includegraphics[width=0.9\linewidth]{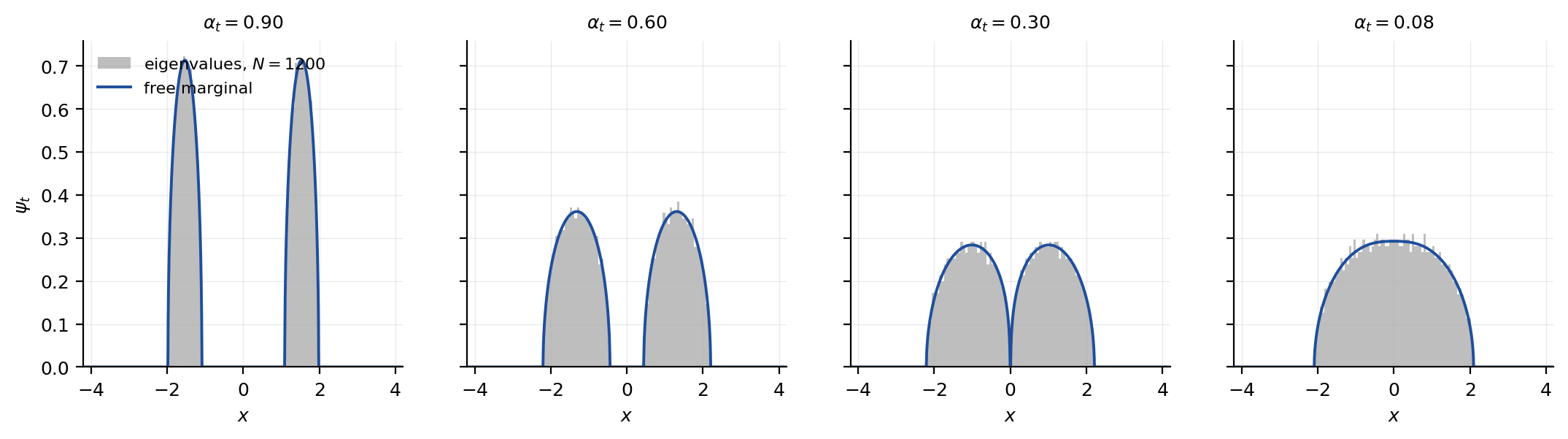}
\caption{Forward free OU (Dyson) flow from the two-atom law. Histograms: a single
Hermitian-diffusion realization, $N=1200$. Curves: the free Fokker--Planck
prediction. The spectrum passes bimodal to unimodal, converging to semicircular.}
\label{fig:forward}
\end{figure}

\begin{figure}[h]
\centering
\includegraphics[width=0.9\linewidth]{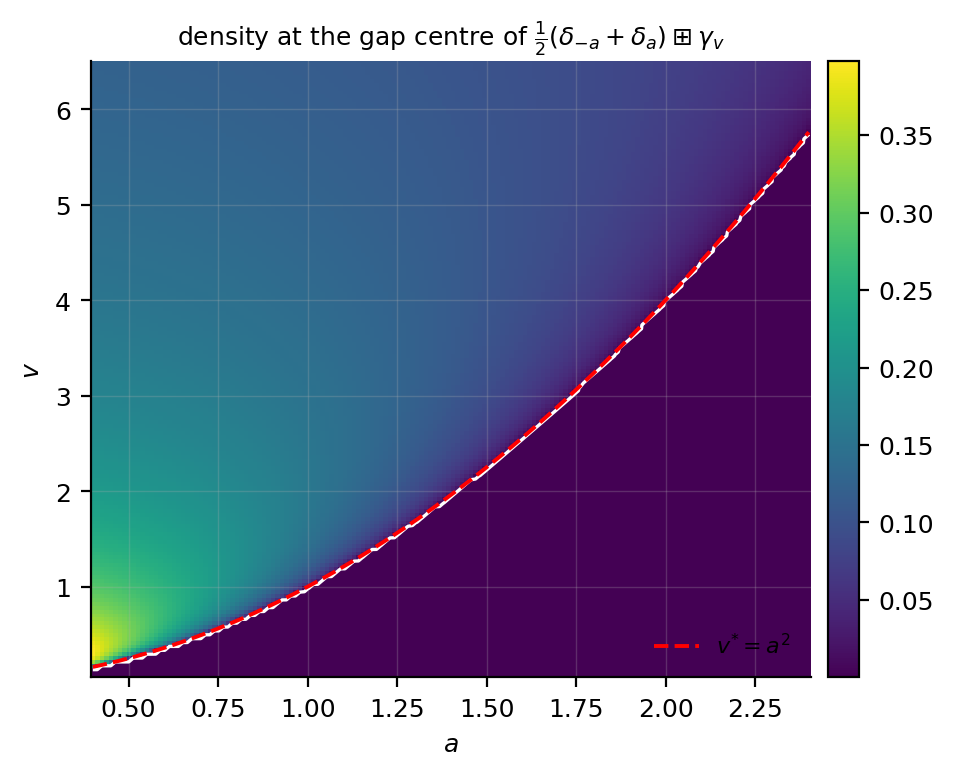}
\caption{Phase diagram of the two-atom support transition. Color: density at the gap
center of $\tfrac12(\delta_{-a}+\delta_a)\boxplus\gamma_v$. White: numerically
detected boundary. Red dashed: the exact transition $v^*=a^2$ of
\S\ref{sec:transition-para}, confirmed as an identity in $a$, not a single point.}
\label{fig:phase}
\end{figure}

\begin{figure}[h]
\centering
\includegraphics[width=0.9\linewidth]{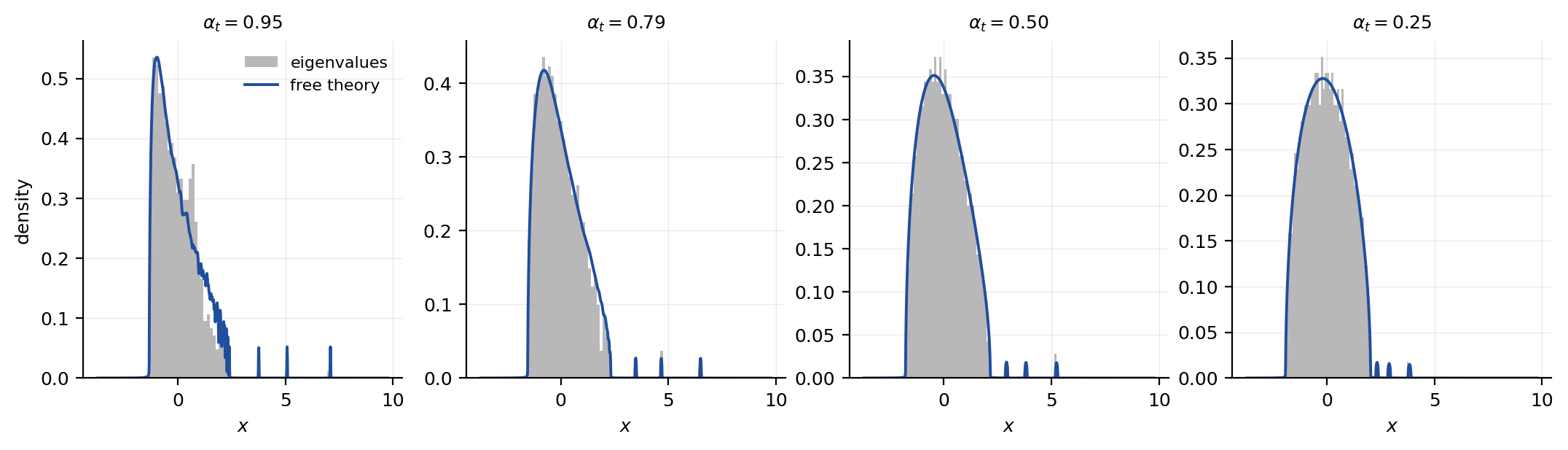}
\caption{Forward free OU flow from a spiked covariance spectrum, $N=600$. The three
outliers merge into the bulk as the flow proceeds -- the finite-rank phenomenon that
Theorem~\ref{thm:gapclosing} does not cover (Appendix~\ref{app:gapclosing}).}
\label{fig:spikedforward}
\end{figure}

\begin{figure}[h]
\centering
\includegraphics[width=0.9\linewidth]{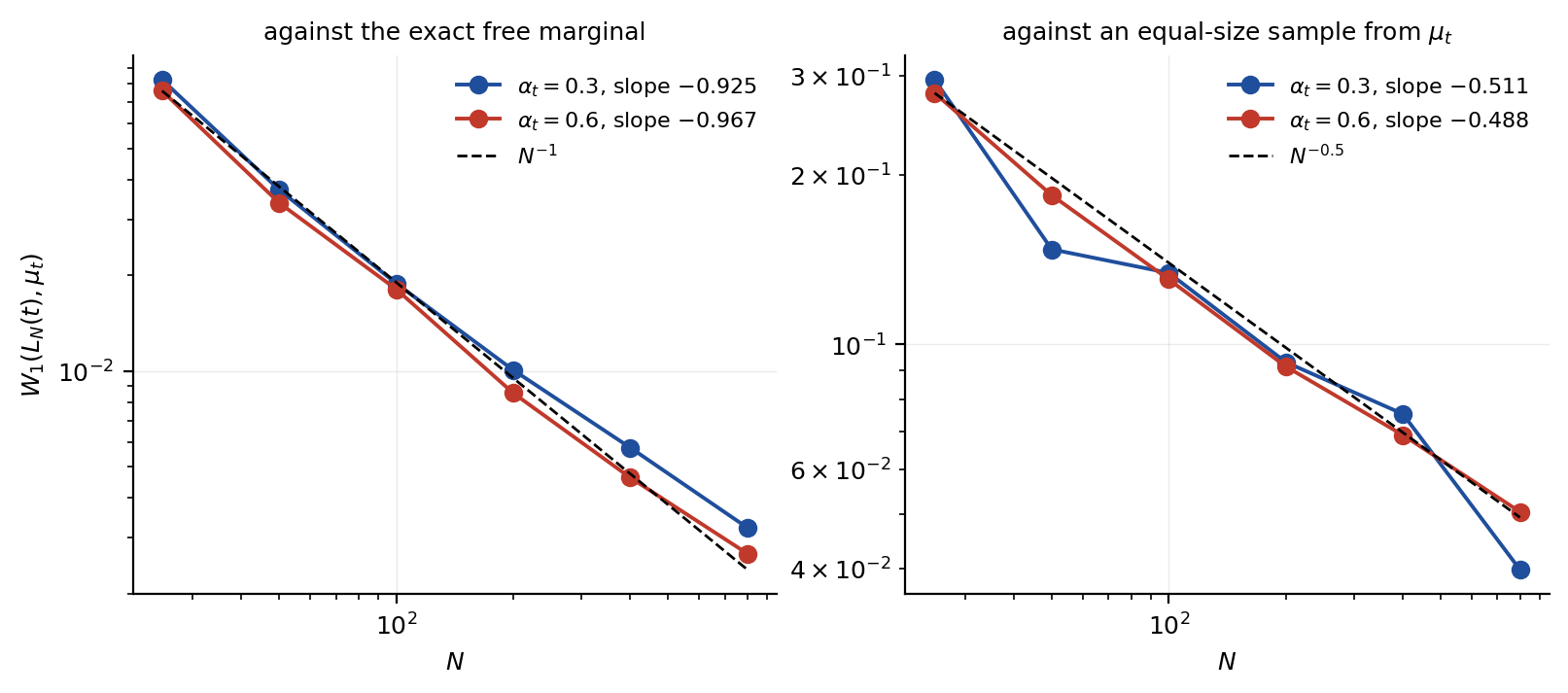}
\caption{Convergence to the free limit, log-log, averaged over independent
realizations. Left: $W_1$ against the exact free marginal, $N^{-1}$ reference. Right:
the same comparison against a finite Monte~Carlo sample instead, $N^{-1/2}$
reference -- the apparent rate is that of the reference sample, not of the matrix
model; see \S\ref{sec:experiments}.}
\label{fig:rate}
\end{figure}

\section{Real-data experiment: full methodology}\label{app:realdata}

\paragraph{Data and preprocessing.} We use \citet{plotlydatasets}'s daily OHLCV data
for the S\&P~500 constituents. Of 505 tickers, 470 have complete coverage over
2013-02-08 to 2018-02-07 (1258 trading days); we keep exactly these, applying no
imputation. Log-returns are computed and each asset's return series is standardized
to unit variance before forming the Pearson correlation matrix on $N=100$ assets (the
first 100 alphabetically by ticker, fixed before looking at any result, to avoid
selecting a flattering subset).

\paragraph{Why the target needs regularizing, and how.} A raw Gaussian kernel density
estimate of the 100 empirical eigenvalues has support all of $\R$, violating the
compact-support hypothesis of Theorem~\ref{thm:design}: precisely the failure mode of
the Cauchy-law example of Appendix~\ref{app:design}'s necessity-of-hypotheses
discussion. We instead truncate the KDE
(bandwidth $0.22$) to $[\min(\mathrm{eigs})-0.3,\max(\mathrm{eigs})+0.3]$ and
renormalize, giving a genuine compactly supported density $\hat\mu_*$, then compute
$b_*=-H\hat\mu_*$ on this interval by direct numerical principal-value integration.
For simulation (not for the stationarity claim itself, which only concerns the
interior of the support) we extend $b_*$ linearly outside the fitted interval with
slope $-5$, a standard confining regularization with no effect on the theory: it only
keeps simulated eigenvalues that fluctuate outside the data-supported range from
escaping instead of returning, exactly as the Marchenko--Pastur experiment
(\S\ref{sec:experiments}) needed a floor near $x=0$ for the same reason.

\paragraph{Dynamical relaxation: what we found, and its limits.} Starting a $N=100$
Hermitian matrix diffusion from an independent Wishart-type spectrum (mismatched in
both shape and scale to the target) and integrating with this drift for $8000$ steps
at $\dd t=5\times10^{-4}$, applying $b_*$ spectrally at each step as in
Algorithm~\ref{alg:main}, the simulated spectrum's bulk moves from its Wishart start
toward the real target's bulk region, but $8000$ steps is not enough for any
eigenvalue to separate into the isolated market-mode outlier at $32.0$: the largest
simulated eigenvalue reaches only $6.4$. This is not a failure of the stationarity
claim -- $\hat\mu_*$, market-mode outlier included, is an exact fixed point of the
continuity equation by construction, exactly as for Example~\ref{ex:mp} -- but a
genuine numerical limitation of the \emph{dynamics}: separating a single eigenvalue
from a bulk of $99$ others against Dyson repulsion is slow, on a timescale this
experiment's budget did not reach. We report this rather than run a longer, possibly
still-inconclusive simulation and describe it as converged: the honest summary is
that the target is provably an equilibrium, the bulk shape is dynamically reached in
the time available, and the isolated outlier's relaxation time is future numerical
work, not a gap in the theory. This mirrors, and sharpens with real rather than
synthetic data, the same honest-reporting standard the Marchenko--Pastur experiment
follows in \S\ref{sec:experiments}: report what is verified, not a single
terminal-distance number that would misrepresent what has actually been checked.

\paragraph{The semicircular baseline.} For contrast, the best-fitting semicircular
law (mean and variance matched to the real spectrum -- the most favorable comparison
possible) misses the true spectrum at $W_1=2.20$, an order of magnitude worse than
any of the synthetic-data $W_1$ values reported in \S\ref{sec:experiments}, and this
mismatch cannot be improved by more training, more data, or a better-tuned schedule:
it is forced by the explicit forward solution of \S\ref{sec:background}, i.e.\ by the
fact that a constant-coefficient free diffusion has no equilibrium other than
semicircular. This
is the concrete, real-data version of the restriction Theorem~\ref{thm:design}
removes.

\section{Limitations and further directions}\label{app:limitations}

\textbf{Edge regularity.} Theorem~\ref{thm:design} requires H\"older continuity of
the target density, and the extension to targets with an $L^2$ but unbounded edge
singularity (as arise for potentials with a logarithmic term at a finite endpoint,
e.g.\ Jacobi-type ensembles) is open.

\textbf{Well-posedness of the designed SDE.} Theorem~\ref{thm:design} establishes
stationarity of the continuity equation; well-posedness of the SDE itself for the
(generally unbounded, non-globally-Lipschitz) designed drift $b_*$ needs the
Lyapunov-type dissipativity conditions of \citet{WeiYin2026} rather than the bounded-
coefficient theory of \citet{BianeSpeicher1998} used elsewhere in this paper; we do
not carry this out (see the discussion following Theorem~\ref{thm:opvol} in the main
text).

\textbf{Multivariate extension.} The results here are for a single self-adjoint
variable; the free transport metric is known to coincide with the classical one on
$\R$ only in this one-variable setting \citep{BianeVoiculescu2001}, and free monotone
transport in several noncommuting variables \citep{GuionnetShlyakhtenko2014} is known
only perturbatively, near a free semicircular family. Whether an analogue of
Theorem~\ref{thm:design} exists for several noncommuting variables is open.

\textbf{Hard edges.} Proposition~\ref{prop:designunique}'s uniqueness statement
requires the target density to vanish at the edges of its support; we do not know an
example where uniqueness fails at a hard edge, only that the proof does not cover it.

\section{Figures moved from the experiments section}\label{app:movedfigs}

\begin{figure}[h]
\centering
\includegraphics[width=0.85\linewidth]{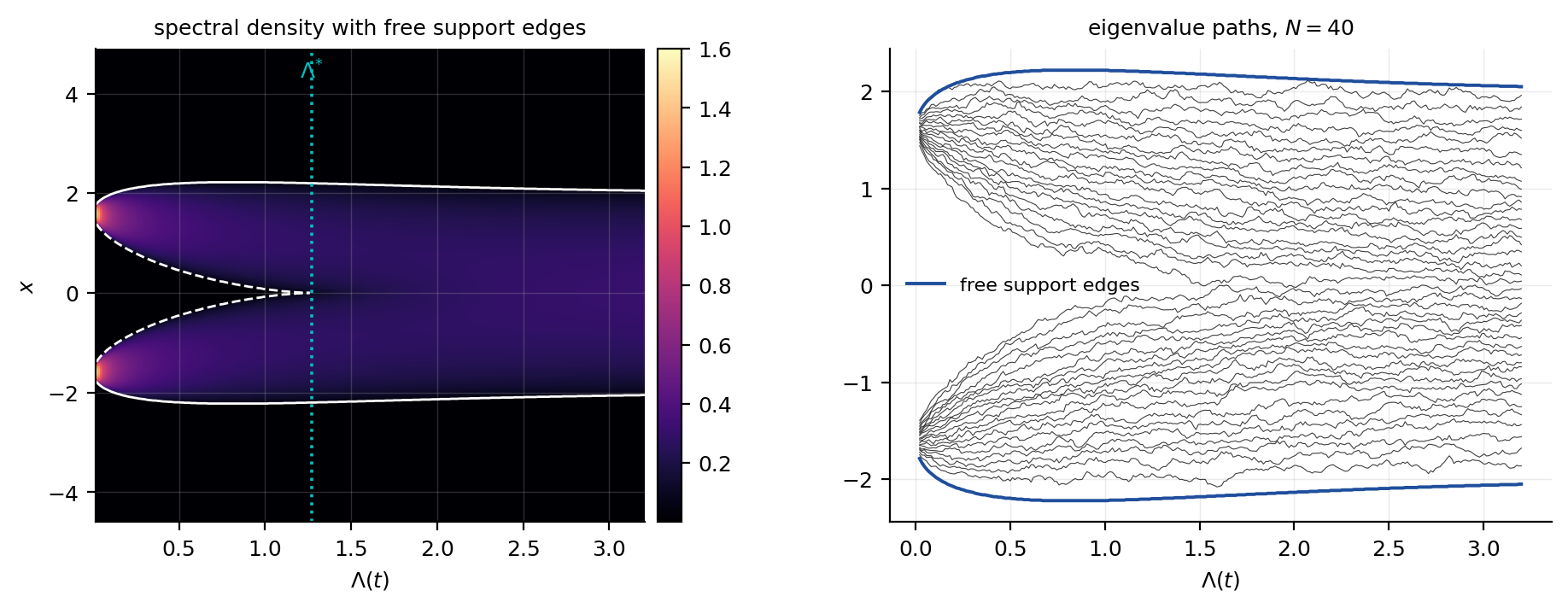}
\caption{Forward flow in spacetime. Left: density $\psi_t(x)$ over $(\Lambda(t),x)$
with support/gap edges from the free theory and the exact transition time
$\Lambda^*$. Right: an $N=40$ eigenvalue trajectory against the same edges.}
\label{fig:spacetime}
\end{figure}

\begin{figure}[h]
\centering
\includegraphics[width=0.85\linewidth]{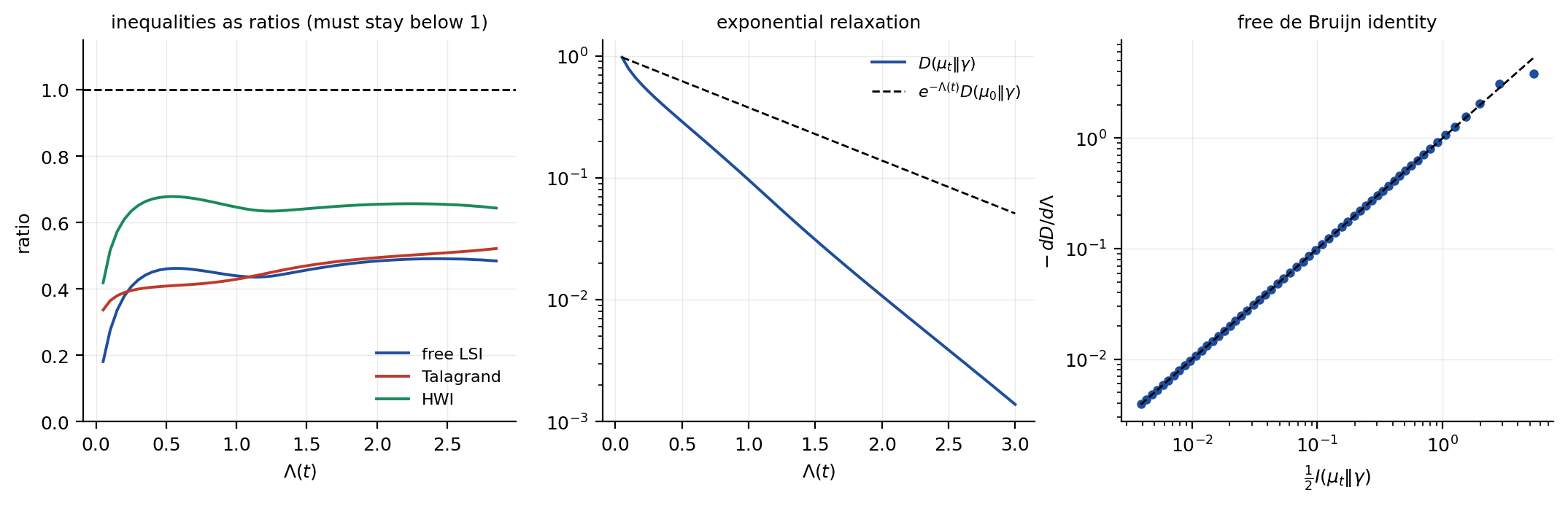}
\caption{Left: the LSI, Talagrand and HWI ratios, all bounded by $1$. Middle:
exponential decay of $\Div$ against $e^{-\Lambda(t)}$. Right: pointwise de~Bruijn
identity $-\dd\Div/\dd\Lambda=\tfrac12\relF$.}
\label{fig:lsi}
\end{figure}

\end{document}